\documentstyle[11pt,amssymb,graphicx]{article}

\begin{document}
\renewcommand{\thefootnote}{}
\date{}

\def\thebibliography#1{\noindent{\normalsize\bf References}
 \list{{\bf
 \arabic{enumi}}.}{\settowidth\labelwidth{[#1]}\leftmargin\labelwidth
 \advance\leftmargin\labelsep
 \usecounter{enumi}}
 \def\newblock{\hskip .11em plus .33em minus .07em}
 \sloppy\clubpenalty4000\widowpenalty4000
 \sfcode`\.=1000\relax}

\baselineskip=12pt
\title{\vspace*{-1cm}
{\large BAND DESCRIPTION OF KNOTS AND VASSILIEV INVARIANTS}}

\author{{\normalsize KOUKI TANIYAMA}\\
{\small Department of Mathematics, College of Arts and Sciences}\\[-1.5mm]
{\small Tokyo Woman's Christian University}\\[-1.5mm]
{\small Zempukuji 2-6-1, Suginamiku, Tokyo 167-8585, Japan}\\[-1.5mm]
{\small e-mail: taniyama@twcu.ac.jp}\\[3mm]
{\normalsize AKIRA YASUHARA }\\
{\small Department of Mathematics, Tokyo Gakugei University}\\[-1.5mm]
{\small Nukuikita 4-1-1, Koganei, Tokyo 184-8501, Japan}\\
{\small {\em Current address}: }\\[-1.5mm]
{\small Department of Mathematics, The George Washington University}\\[-1.5mm]
{\small Washington, DC 20052, USA}\\[-1.5mm]
{\small e-mail: yasuhara@u-gakugei.ac.jp}\\
}

\maketitle

\vspace*{-5mm}  

{\small 
\begin{quote}
\begin{center}A{\sc bstract}\end{center}
\baselineskip=10pt
\hspace*{1em}In 1993 K. Habiro defined {$C_k$-move}
of oriented links and around
1994 he proved that two oriented knots are transformed into 
each other by $C_k$-moves if and only if
they have the same Vassiliev invariants of order $\leq k-1$. 
In this paper we define Vassiliev invariant of type $(k_1,...,k_l)$, 
and show that, for $k=k_1+\cdots+k_l$, two oriented knots are 
transformed into each other by $C_k$-moves if and only if
they have the same Vassiliev invariants of type $(k_1,...,k_l)$. 
We introduce a concept \lq band description of knots' 
and give a diagram-oriented proof of this theorem. 
When $k_1=\cdots=k_l=1$, the Vassiliev invariant of type
$(k_1,...,k_l)$ coincides
with the Vassiliev invariant of order $\leq l-1$ in the usual sense. 
As a special case, we have Habiro's theorem stated above. 
\end{quote}}

\baselineskip=12pt


\noindent
{\bf Introduction}

\medskip

In 1993 K. Habiro defined {\it $C_k$-move} of oriented links for each
natural number $k$ \cite{Habiro2}. A $C_k$-move is a kind of local move 
of oriented links. Around 1994 he proved that two oriented
knots have the same Vassiliev invariants of order 
$\leq k-1$ if and only if they are transformed into each other by
$C_k$-moves. Thus he has succeeded in deducing a geometric conclusion from
an algebraic
condition. However this theorem appears only in his recent paper
\cite{Habiro1}. In \cite{Habiro1} he develops his original 
clasper theory and obtains the theorem 
as a consequence of clasper theory.  We note that the \lq
if' part of the theorem is also shown in
\cite{Gusarov2}, \cite{Ohyama}, \cite{Stanford0} and \cite{T-Y}, and 
in \cite{Stanford3} T. Stanford gives another characterization of 
knots with the same Vassiliev invariants of order  $\leq k-1$.

In this paper we define Vassiliev invariant of type $(k_1,...,k_l)$, 
and show that, for $k=k_1+\cdots+k_l$, two oriented knots are 
transformed into each other by $C_k$-moves if and only if
they have the same Vassiliev invariants of type 
$(k_1,...,k_l)$. When $k_1=\cdots=k_l=1$, the Vassiliev invariant 
of type $(k_1,...,k_l)$ coincides 
with the Vassiliev invariant of order $\leq l-1$ in the usual sense. 
As a special case, 
we have Habiro's theorem. We use 
a concept \lq band description of knots' 
and give a diagram-oriented proof of the theorem. 
The proof is elementary and completely self-contained. 
Note that the prototypes of band description appear in
\cite{Suzuki}, \cite{Yamamoto} and \cite{Yasuhara}. In particular in
\cite{Yasuhara} the second author showed that any knot can be expressed 
as a band sum of a trivial knot and some Borromean rings. The
concept of band description is a development of this fact. 
More generally the authors defined \lq band
description of spatial graphs' in \cite{T-Y}. 
The concept 
of Vassiliev invariant of type $(k_1,...,k_l)$ is defined in \cite{T-Y}. 
A related result to the case $k_1=\cdots=k_l=2$ is shown in
\cite{Stanford4}.

\bigskip\noindent
{\bf 1. Definitions and Main Result}

\medskip

Throughout this paper we work in the piecewise linear category.

A {\it tangle} $T$ is a disjoint union of properly embedded 
arcs in the unit $3$-ball $B^{3}$. A tangle $T$ is {\it trivial} if there 
exists a properly embedded disk in $B^3$ that contains $T$. 
A {\it local move} is a pair of trivial tangles 
$(T_{1},T_{2})$ with $\partial T_{1}=\partial T_{2}$ 
such that for each component $t$ of $T_1$ there exists 
a component $u$ of $T_2$ with $\partial t=\partial u$. Such a pair of
components is called a {\it
corresponding pair}. Two local moves $(T_{1},T_{2})$ and $(U_{1},U_{2})$
are {\it equivalent}, 
denoted by $(T_{1},T_{2})\cong (U_{1},U_{2})$, 
if there is an orientation preserving 
self-homeomorphism $\psi :B^{3}\rightarrow B^{3}$ such that $\psi (T_{i})$ 
and $U_{i}$ are ambient isotopic in $B^3$ relative to $\partial
B^{3}$ for 
$i=1,2$. Here $\psi (T_{i})$ 
and $U_{i}$ are {\it ambient isotopic in $B^3$ relative to $\partial
B^{3}$} if $\psi (T_{i})$ is deformed to $U_{i}$ by an ambient isotopy of
$B^3$ that is
pointwisely fixed on $\partial B^3$.

Let $(T_{1},T_{2})$ be a 
local move, $t_{1}$ a component of $T_1$ and $t_2$ a component of $T_2$ 
with $\partial t_{1}=\partial t_{2}$. Let $N_1$ and $N_2$ be regular 
neighbourhoods of $t_1$ and $t_2$ in $(B^3-T_1)\cup t_1$ and $(B^3-T_2)\cup
t_2$ respectively such
that $N_{1}\cap
\partial  B^{3}=N_{2}\cap \partial B^{3}$. Let $\alpha$ be a disjoint union
of properly 
embedded arcs in $B^{2}\times [0,1]$ as illustrated in Fig. 1.1.
Let $\psi_{i}:B^{2}\times [0,1]\rightarrow N_{i}$ be a homeomorphism 
with $\psi_{i}(B^{2}\times \{ 0,1\} )=N_{i}\cap \partial B^{3}$ for $i=1,2$. 
Suppose that $\psi_{1}(\partial \alpha )=\psi_{2}(\partial \alpha )$ and 
$\psi_{1}(\alpha )$ and $\psi_{2}(\alpha )$ are ambient isotopic in $B^{3}$ 
relative to $\partial B^3$. Then we say that a local move 
$((T_{1}-t_{1})\cup \psi_{1}(\alpha ), (T_{2}-t_{2})\cup \psi_{2}(\alpha ))$ 
is a {\it double of $(T_{1},T_{2})$ with respect to the components $t_1$ 
and $t_2$}. 
Note that a double of $(T_{1},T_{2})$ with respect to $t_1$ 
and $t_2$ is well-defined up to equivalence.

\begin{center} 
\includegraphics[trim=0mm 0mm 0mm 0mm, width=.2\linewidth]
{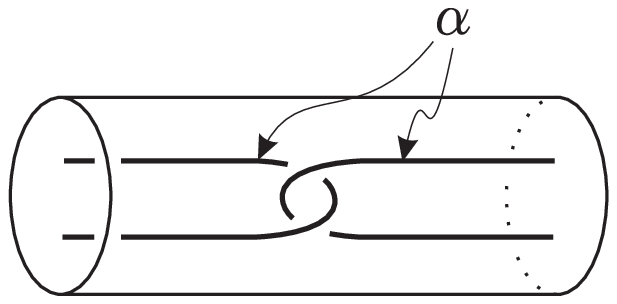}

Fig. 1.1
\end{center}

A {\it $C_{1}$-move} is a local move $(T_1,T_2)$ as illustrated in Fig. 1.2. 
A double of a $C_{k}$-move is called a {\it $C_{k+1}$-move}. 
Note that for each natural number $k$ there are only finitely 
many $C_{k}$-moves up to 
equivalence. We note that by the definition a $C_k$-move is {\it Brunnian}.
That is, if $(T_1,T_2)$ is a 
$C_k$-move and $t_1$, $t_2$ components of $T_1$ and $T_2$ respectively with
$\partial t_1=\partial
t_2$, then the tangles $T_1-t_1$ and $T_2-t_2$ are ambient isotopic in
$B^3$ relative to $\partial B^3$.

\begin{center}
\includegraphics[trim=0mm 0mm 0mm 0mm, width=.3\linewidth]
{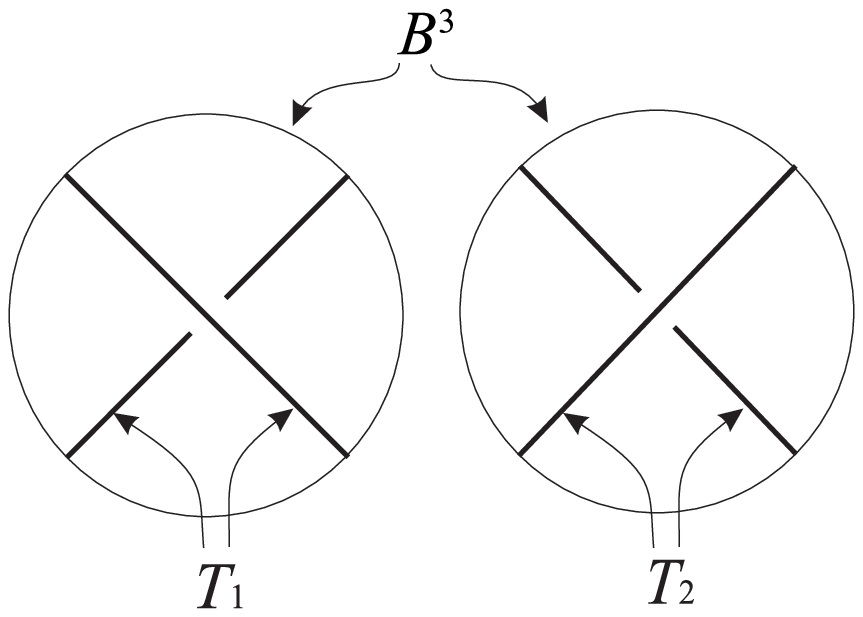}

Fig. 1.2
\end{center}

Let $(T_{1},T_{2})$ be a local move. Then $(T_{2},T_{1})$ is also a local 
move. We call $(T_{2},T_{1})$ the {\it inverse} of $(T_{1},T_{2})$. It is 
easy to see that the inverse of a $C_{1}$-move is equivalent to itself. 
Then it follows inductively that the inverse of a $C_{n}$-move is equivalent 
to a $C_{n}$-move (but possibly not equivalent to itself).

Let $K_1$ and $K_2$ be oriented knots in the oriented three-sphere $S^3$. 
We say that $K_1$ and
$K_2$ are  {\it related by a local move $(T_{1},T_{2})$} if 
there is an orientation preserving embedding 
$h:B^{3}\rightarrow S^{3}$ such that $K_{i}\cap h(B^{3})=h(T_{i})$ for 
$i=1,2$ and $K_{1}-h(B^{3})=K_{2}-h(B^{3})$ together with orientations. 
Then we also say that $K_2$ 
is obtained from $K_1$ by an {\it application of $(T_{1},T_{2})$}. Two
oriented knots $K_{1}$ and
$K_{2}$ are {\it $C_{k}$-equivalent} if 
$K_{1}$ and $K_{2}$ are related by a finite sequence of $C_{k}$-moves 
and ambient isotopies. This relation is an equivalence relation on knots. 
It is known that $C_k$-equivalence implies 
$C_{k-1}$-equivalence \cite{Habiro1}, see Remark 2.4.

Let $l$ be a positive integer and $k_1,...,k_l$ positive integers. Suppose
that for each
$P\subset\{1,...,l\}$ an oriented knot $K_P$ in $S^3$ is assigned. Suppose 
that there are orientation preserving embeddings 
$h_i:B^3\rightarrow S^3$ $(i=1,...,l)$ such that \\
(1) $h_i(B^3)\cap h_j(B^3)=\emptyset$ if $i\neq j$,\\
(2) $K_P-\bigcup_{i=1}^l h_i(B^3)=K_{P'}-\bigcup_{i=1}^l h_i(B^3)$ together
with orientation 
for any subsets $P,P'\subset \{1,...,l\}$,\\
(3) $(h_i^{-1}(K_\emptyset),h_i^{-1}(K_{\{1,...,l\}}))$ is 
a $C_{k_i}$-move $(i=1,...,l)$, and\\
(4) $K_P\cap h_i(B^3)=\left\{
\begin{array}{ll}
K_{\{1,...,l\}}\cap h_i(B^3) & \mbox{if $i\in P$},\\
K_\emptyset\cap h_i(B^3) & \mbox{otherwise}.
\end{array}
\right.$\\
Then we call the set of knots $\{K_P|P\subset\{1,...,l\}\}$ a {\em singular
knot of type
$(k_1,...,k_l)$}. Let ${\cal K}$ be the set of all oriented knot types in
$S^3$ and ${\Bbb Z}{\cal K}$
the free abelian group  generated by ${\cal K}$. We sometimes identify a
knot and its knot type
without explicite mention. For a singular knot
$K=\{K_P|P\subset\{1,...,l\}\}$ of type
$(k_1,...,k_l)$, we define an element $\kappa(K)$ of ${\Bbb Z}{\cal K}$ by 
\[\kappa(K)=\sum_{P\subset\{1,...,l\}}(-1)^{| P|}K_P.\]
Let ${\cal V}(k_1,...,k_l)$ be the subgroup of ${\Bbb Z}{\cal K}$
generated by all $\kappa(K)$ where $K$ varies over all
singular knots of type $(k_1,...,k_l)$.

Let $K_1\#K_2$ be the composite knot of two knots $K_1$ and $K_2$. Then
$K_1\#K_2-K_1-K_2\in{\Bbb
Z}{\cal K}$ is called a {\em composite relator} following Stanford
\cite{Stanford3}. Let ${\cal R}_\#$
be the subgroup of
${\Bbb Z}{\cal K}$  generated by all composite relators.

Let $\iota:{\cal K}\rightarrow{\Bbb Z}{\cal K}$ be the natural 
inclusion map.
Let $\pi_{(k_1,...,k_l)}:{\Bbb Z}{\cal K}\rightarrow
{\Bbb Z}{\cal K}/{\cal V}(k_1,...,k_l)$ and $\lambda_{(k_1,...,k_l)}:{\Bbb
Z}{\cal K}/{\cal
V}(k_1,...,k_l)\rightarrow{\Bbb Z}{\cal K}/({\cal
V}(k_1,...,k_l)+{\cal R}_\#)$ be the quotient homomorphisms. Then the
composite maps
$\pi_{(k_1,...,k_l)}\circ\iota:
{\cal K}\rightarrow {\Bbb Z}{\cal K}/{\cal V}(k_1,...,k_l)$ and 
$\lambda_{(k_1,...,k_l)}\circ\pi_{(k_1,...,k_l)}\circ\iota:
{\cal K}\rightarrow {\Bbb Z}{\cal K}/({\cal V}(k_1,...,k_l)+{\cal R}_\#)$ 
are called the {\em universal Vassiliev invariant of type 
$(k_1,...,k_l)$} and the {\em universal additive Vassiliev invariant of type 
$(k_1,...,k_l)$} respectively. We denote them by $v_{(k_1,...,k_l)}$ and
$w_{(k_1,...,k_l)}$
respectively.

Since a $C_1$-move is a crossing change we have that a singlar knot of type
$(\underbrace{1,...,1}_{l})$ is essentially the same as a singlar knot with
$l$ crossing vertices
in the usual sense. Therefore we have that
$\displaystyle  v_{(1,...,1)}$ is the universal 
Vassiliev invariant of order $\leq l-1$ and $\displaystyle 
w_{(1,...,1)}$ is the universal additive
Vassiliev invariant of order $\leq l-1$. Note that $\displaystyle 
v_{(1,...,1)}(K_1)=\displaystyle 
v_{(1,...,1)}(K_2)$ if and only if $v(K_1)=v(K_2)$ for any Vassiliev
invariant $v$ of order $\leq l-1$. Similarly
$w_{(1,...,1)}(K_1)=\displaystyle 
w_{(1,...,1)}(K_2)$ if and only if $w(K_1)=w(K_2)$ for any
additive Vassiliev invariant $w$ of order $\leq l-1$.
We also note that $v_{(2,...,2)}$ is essentially same as that 
defined in \cite{Mellor}, \cite{Stanford4}. 
In \cite{T-Y} the authors 
defined a finite type invariant of order $(k;n)$, 
that is essentially same as 
$\displaystyle v_{(\scriptsize\underbrace{n-1,...,n-1}_{k+1})}$.

Now we state our main results.

\medskip
\noindent
{\bf Theorem 1.1.} {\em Let $k_1,...,k_l$ be positive integers 
and $k=k_1+\cdots+k_l$. Then 
${\cal V}(k)\subset {\cal V}(k_1,...,k_l)$.}

\medskip
\noindent
{\bf Theorem 1.2.} {\em Let $k_1,...,k_l$ be positive integers 
and $k=k_1+\cdots+k_l$. Then ${\cal V}(k_1,...,k_l)\subset 
{\cal V}(k)+{\cal R}_\#$.}

By Theorems 1.1 and 1.2, we have the following corollary.

\medskip
\noindent
{\bf Corollary 1.3.} {\em Let $k_1,...,k_l$ be positive integers 
and $k=k_1+\cdots+k_l$. Then ${\cal V}(k_1,...,k_l)+{\cal R}_\#=
{\cal V}(k)+{\cal R}_\#$. $\Box$}

\medskip
The following theorem was proved by Habiro \cite{Habiro1}, and 
the authors gave a 
proof by using band description \cite{T-Y}. In section 3, we give the 
same proof as in \cite{T-Y} for the convenience of the reader.

\medskip
\noindent
{\bf Theorem 1.4.} (Habiro \cite{Habiro1}) {\it The
$C_k$-equivalence classes of oriented knots in
$S^3$ forms an abelian group under connected sum of oriented knots.}

\medskip
We denote this group by ${\cal K}/C_k$.

\medskip
\noindent
{\bf Theorem 1.5.} {\it Let $\eta_k:{\cal K}/C_k\rightarrow{\Bbb Z}{\cal
K}/({\cal
V}(k)+{\cal R}_\#)$ be a map induced by the inclusion $\iota$. Then
$\eta_k$ is an isomorphism.}

\medskip
Let $K_1$ and $K_2$ be oriented knots and $k=k_1+\cdots+k_l$. 
If $K_1$ and $K_2$ are $C_k$-equivalent, then by Theorem 1.1, 
$K_1-K_2\in{\cal V}(k)\subset {\cal V}(k_1,...,k_l)$. 
Therefore we have $v_{(k_1,...,k_l)}(K_1)=v_{(k_1,...,k_l)}(K_2)$. 
On the other hand, if
$v_{(k_1,...,k_l)}(K_1)=v_{(k_1,...,k_l)}(K_2)$, then
$w_{(k_1,...,k_l)}(K_1)=w_{(k_1,...,k_l)}(K_2)$. Then
by Corollary 1.3, $K_1-K_2\in{\cal V}(k)+{\cal R}_\#$. Then by Theorem 1.5 we
have that 
$K_1$ and $K_2$ are $C_k$-equivalent. Hence we have the following 
theorem. 

\medskip
\noindent
{\bf Theorem 1.6.} {\em Let $k_1,...,k_l$ be positive integers 
and $k=k_1+\cdots+k_l$. Let $K_1$ and $K_2$ be oriented knots in $S^3$.
Then the following conditions
are mutually equivalent.

{\rm (1)} $K_1$ and $K_2$ are $C_k$-equivalent,

{\rm (2)} $v_{(k_1,...,k_l)}(K_1)=v_{(k_1,...,k_l)}(K_2)$,

{\rm (3)} $w_{(k_1,...,k_l)}(K_1)=w_{(k_1,...,k_l)}(K_2)$.
$\Box$}

\medskip
As a special case of Theorem 1.6 we have the following theorem.

\medskip
\noindent{\bf Theorem 1.7.} (Habiro \cite{Habiro1}) 
{\it Two oriented knots $K_1$ and $K_2$ are $C_{k}$-equivalent
if and only if their values of the universal {\rm(}additive{\rm)} Vassiliev
invariant 
of order $\leq k-1$ are equal.  $\Box$}

\medskip
The remainder of this paper, we prove Theorems 1.1, 1.2, 1.4 and 1.5.
We give a proof of Theorem 1.1 in section 2 and proofs of Theorems 1.2, 
1.4 and 1.5 in section 3. The reader who wishes may change the order 
of sections 2 and 3 since they are independent except for Remark 2.4.

\bigskip\noindent
{\bf 2. $C_k$-moves}

\medskip
A graph is called a {\em tree} if it is connected and simply connected as a
topological space. A
graph is {\it uni-trivalent} if each vertex has degree one or three. Let
$G$ be a uni-trivalent tree
embedded on the unit disk $D^2$ such that $G\cap\partial D^2$ is exactly
the set of degree-one
vertices $\{v_1,...,v_{k+1}\}$ of $G$. Suppose that an edge $e$ of $G$ is
specified. We will assign 
a $C_k$-move to $G$ with respest to $e$ 
and a pair of corresponding components of it to each
$v_i$ as follows. If $G$ is a tree with just two vertices
$v_1,v_2$ and one edge
$e$ joining them then we assign a diagram of a $C_1$-move as illustrated in
Fig. 2.1. Suppose that
for each uni-trivalent tree on $k$ vertices with a specified edge, 
a diagram of a $C_{k-1}$-move is assigned. Suppose that 
$v_k$ and $v_{k+1}$ are degree-one vertices of $G$ such that there is a
degree-three vertex $u$ of
$G$ that is adjacent to both of $v_k$ and $v_{k+1}$. Suppose that 
neither $uv_k$ nor $uv_{k+1}$ is a specified edge $e$. 
Let $G'$ be a uni-trivalent tree obtained
from $G$ by deleting $v_{k+1}$ and the edge $uv_{k+1}$ and
forgetting $u$. Let ${\cal
D'}$ be a diagram of a $C_{k-1}$-move assigned to $G'$ 
with respect to $e$. Let $t_1$ and $t_2$ be a pair of
corresponding components assigned to $v_k$. Then we replace their parts in
$D'$ as the diagram
illustrated in Fig. 2.2. We assign the new pairs of corresponding
components to $v_k$ and $v_{k+1}$
respecting the cyclic order of them on $\partial D^2$. Thus we have
assigned a $C_k$-move to $G$ with respect to the specified edge.
See for example Fig. 2.3. 
Note that any $C_k$-move is assigned to a uni-trivalent tree with respect 
to a specified edge up to equivalence of local moves. 

\begin{center}
\includegraphics[trim=0mm 0mm 0mm 0mm, width=.5\linewidth]
{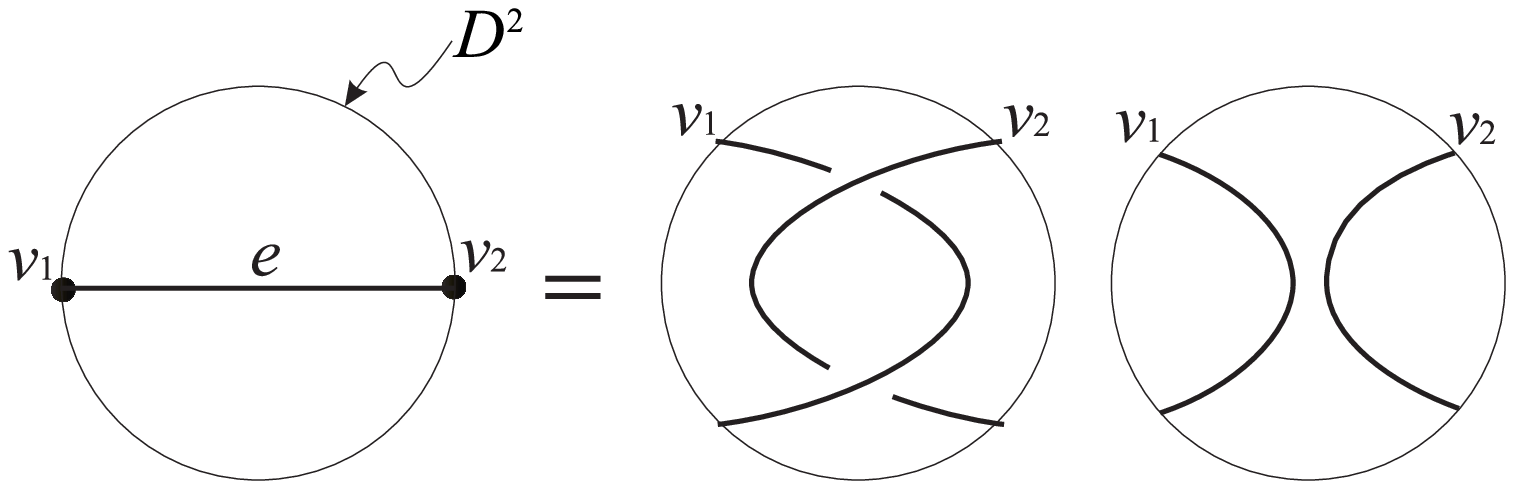}

Fig. 2.1
\end{center}

\begin{center}
\includegraphics[trim=0mm 0mm 0mm 0mm, width=.5\linewidth]
{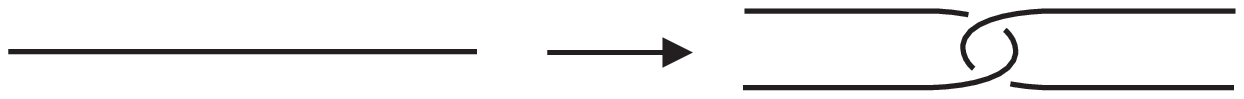}

Fig. 2.2
\end{center}

\begin{center}
\includegraphics[trim=0mm 0mm 0mm 0mm, width=.7\linewidth]
{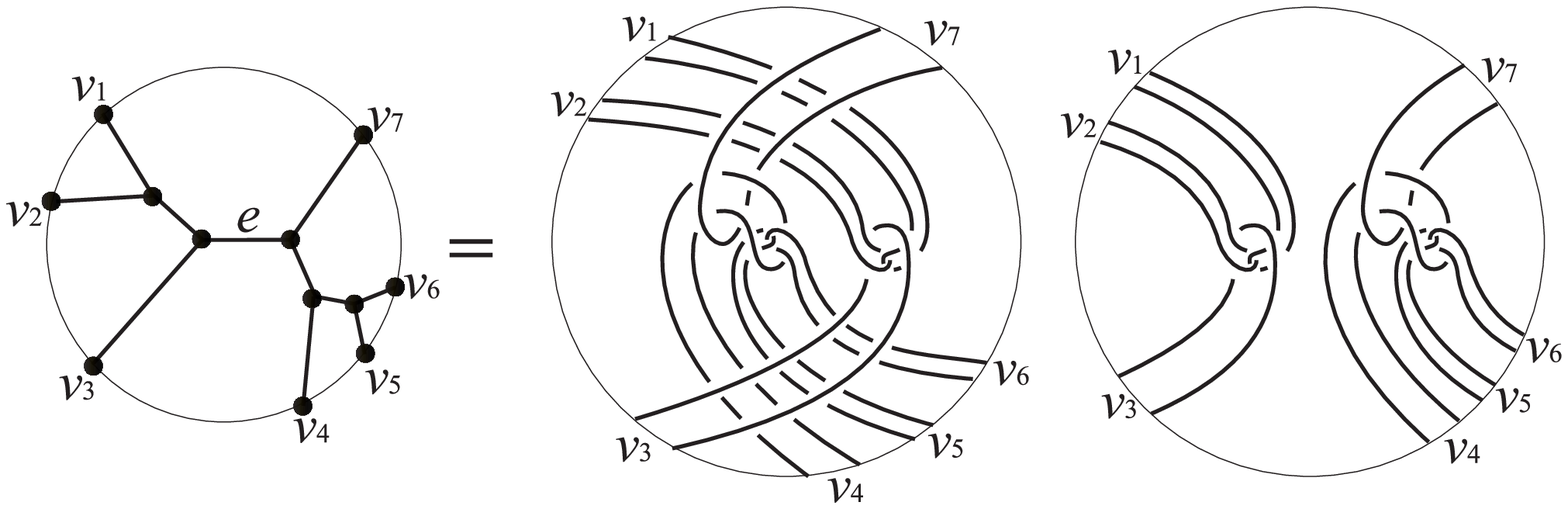}

Fig. 2.3
\end{center}

\medskip
\noindent{\bf Lemma 2.1.} {\it Let $(T_1,T_2)$ be a $C_k$-move assigned to
a uni-trivalent tree $G$ with respect to a specified edge $e$. 
Suppose that an edge $e'$ is
incident to $e$. Then there is
a re-embedding $f:G\rightarrow D^2$ such that the $C_k$-move assigned to
$f(G)$ with respect to
$f(e')$ is equivalent to $(T_1,T_2)$.}

\medskip
\noindent{\bf Proof.} Let $v$ be the common vertex of $e$ and $e'$. First
suppose that only $v$ is
the degree-three vertex of $G$. Then we have the result by the deformation
illustrated in Fig. 2.4.
By taking appropriate doubles we have the general case. $\Box$

\begin{center}
\includegraphics[trim=0mm 0mm 0mm 0mm, width=.65\linewidth]
{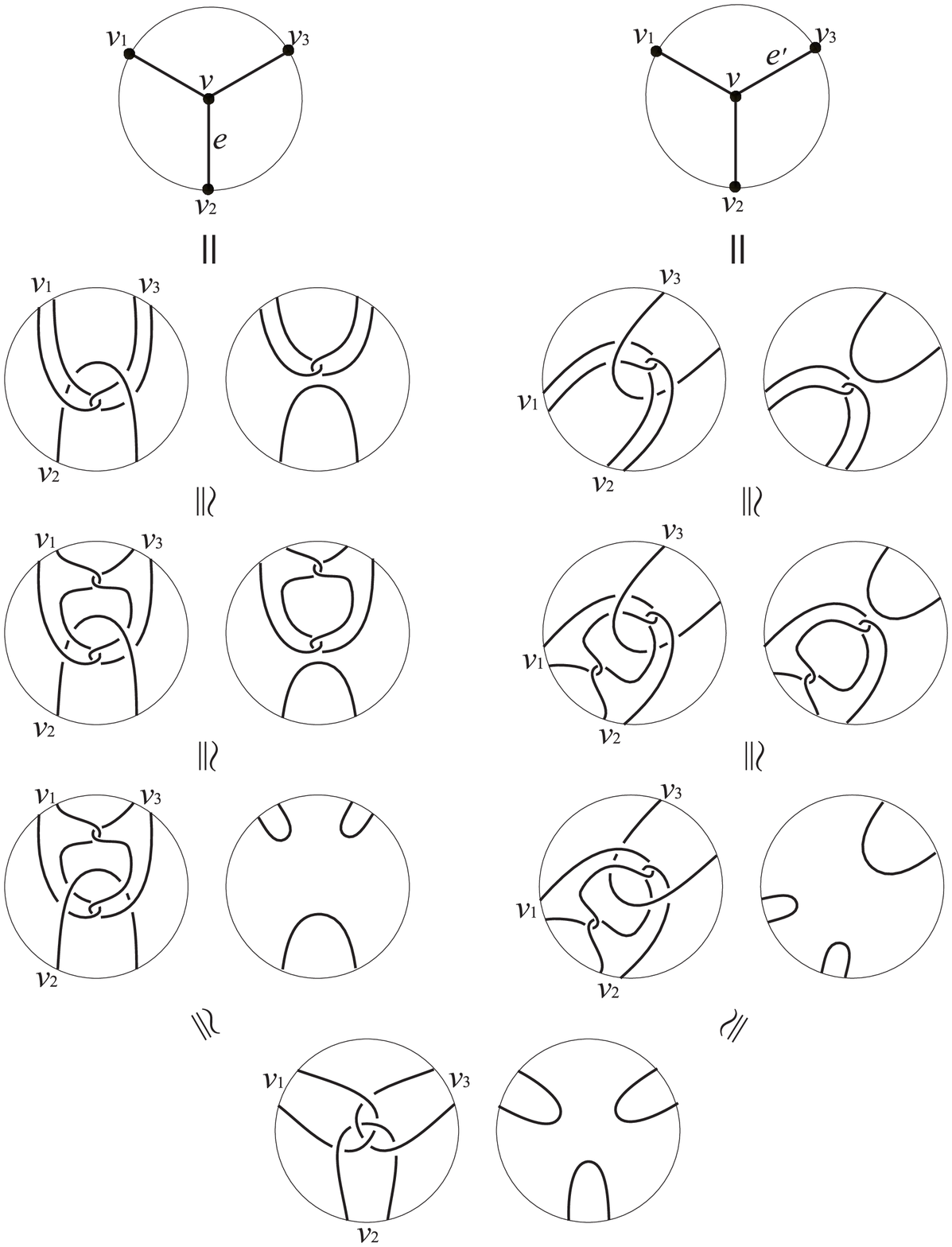}

Fig. 2.4
\end{center}

\medskip
A {\it path} of $G$ is a subgraph of $G$ homeomorphic to a closed interval. Let
$\sigma(G)$ be the
maximal of the number of the edges of a path of $G$ and call it the 
{\it diameter} of $G$. We say that a $C_k$-move
$(T_1,T_2)$ is {\it one-branched} if it is assigned to a uni-trivalent tree
of diameter $k$. 
Note that if $G$ is a uni-trivalent tree with $k+1$ degree-one vertices 
and $\sigma(G)=k$, then all of the degree-three vertices lie on 
the path with $k$ edges.

Let $S_1$ and $S_2$ be tangles.  We say that $S_1$ and
$S_2$ are  {\it related by a local move $(T_{1},T_{2})$} if 
there is an orientation preserving embedding 
$h:B^{3}\rightarrow {\rm int}B^{3}$ such that $S_{i}\cap h(B^{3})=h(T_{i})$
for $i=1,2$ and $S_{1}-h(B^{3})=S_{2}-h(B^{3})$.
Then we say that $S_2$
is obtained from $S_1$ by an {\it application of $(T_{1},T_{2})$}.

The following result was shown by Habiro \cite{Habiro2}. 
Since the article \cite{Habiro2} is written in Japanese, 
we give a proof by using our terms. 

\medskip
\noindent{\bf Lemma 2.2.} (Habiro \cite{Habiro2}) {\it Let $(T_1,T_2)$ be a
$C_k$-move. Then $T_1$
and $T_2$ are related by a finite sequence of one-branched $C_k$-moves and
ambient isotopies relative
to $\partial B^3$.}

\medskip
\noindent{\bf Proof.} Suppose
that $(T_1,T_2)$ is assigned to a uni-trivalent tree
$G$. If
$\sigma(G)=k$ then $(T_1,T_2)$ itself is a one-branched $C_k$-move. 
Therefore we
may suppose that $\sigma(G)<k$. It is sufficient to show that
$T_1$ and $T_2$ are related by
$C_k$-moves each of which is assigned to a uni-trivalent tree of diameter
$\sigma(G)+1$. Let $P$ be
a path of $G$ containing $\sigma(G)$ edges. 
Since $\sigma(G)<k$, there are degree-three vertices of
$G$ that are not on $P$. Let $v$ be one of them such that $v$ is adjacent
to a vertex $w$ on $P$.
Let $w_1$ and $w_2$ be the vertices on $P$ adjacent to $w$. Let $v_1$ and
$v_2$ be the other vertices
adjacent to $v$. By Lemma 2.1 we may assume that the specified edge is
$ww_1$. We temporarily
forget the embedding of $G$ into $D^2$ and define two abstract graphs $G_1$
and $G_2$ from $G$ as
follows. Let $G_i$ be a uni-trivalent tree obtained from $G$ by deleting
the edge
$vv_i$, forgetting $v$, adding a vertex $u$ on $ww_2$ and adding an edge
$uv_i$ $(i=1,2)$. Then we have
$\sigma(G_1)=\sigma(G_2)=\sigma(G)+1$. We will show that $T_1$ and $T_2$
are related by a $C_k$-move
assigned to some embedding of $G_1$ and a $C_k$-move assigned to some
embedding of $G_2$. In the
following we consider the simplest case illustrated in Fig. 2.5.
General cases follow immediately
by taking appropriate doubles. 
The $C_k$-move of Fig. 2.5 is assigned to the graph of Fig. 2.5 with
respect to $ww_1$. Note that the $C_k$-move of Fig. 2.5 is
equivalent to a local
move of Fig. 2.6. Then it is realized by two local moves Fig. 2.7 (a)
and (b). Then each of them is equivalent to a desired $C_k$-move. 
Fig. 2.8 indicates it for the case (a). 
The case (b) is similar and we omit it.
$\Box$

\begin{center}
\includegraphics[trim=0mm 0mm 0mm 0mm, width=.65\linewidth]
{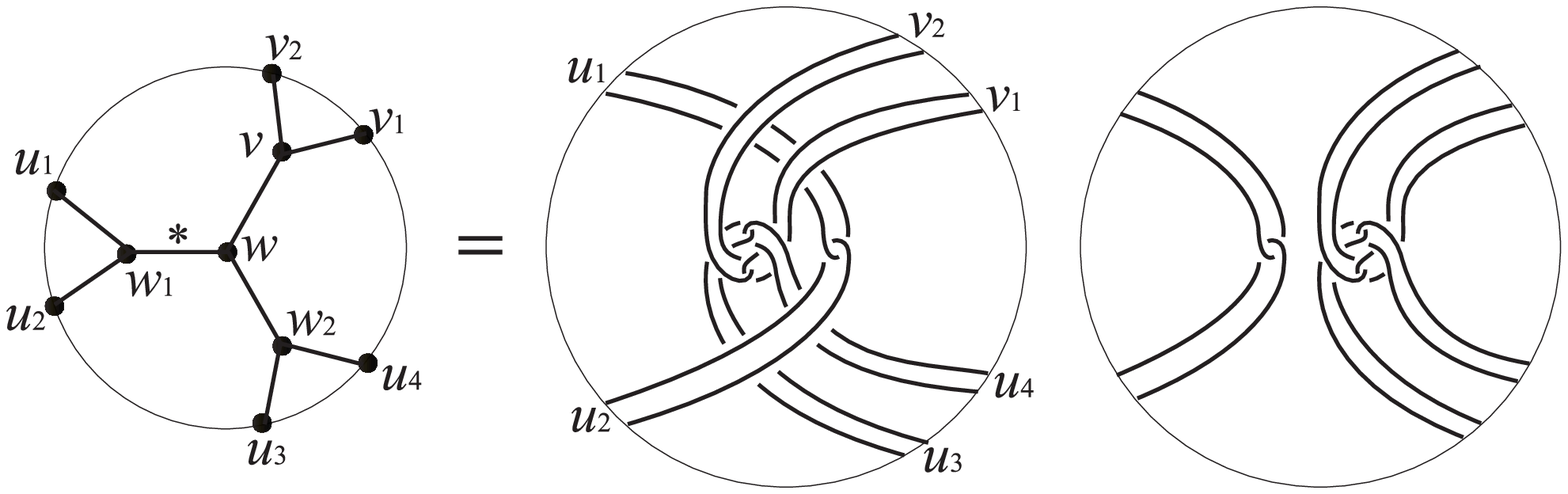}

Fig. 2.5
\end{center}

\begin{center}
\includegraphics[trim=0mm 0mm 0mm 0mm, width=.35\linewidth]
{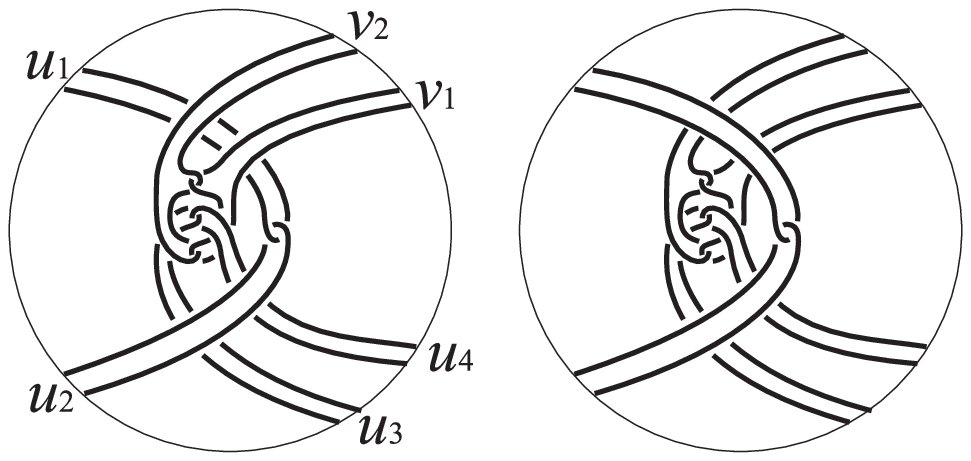}

Fig. 2.6
\end{center}

\begin{center}
\includegraphics[trim=0mm 0mm 0mm 0mm, width=.75\linewidth]
{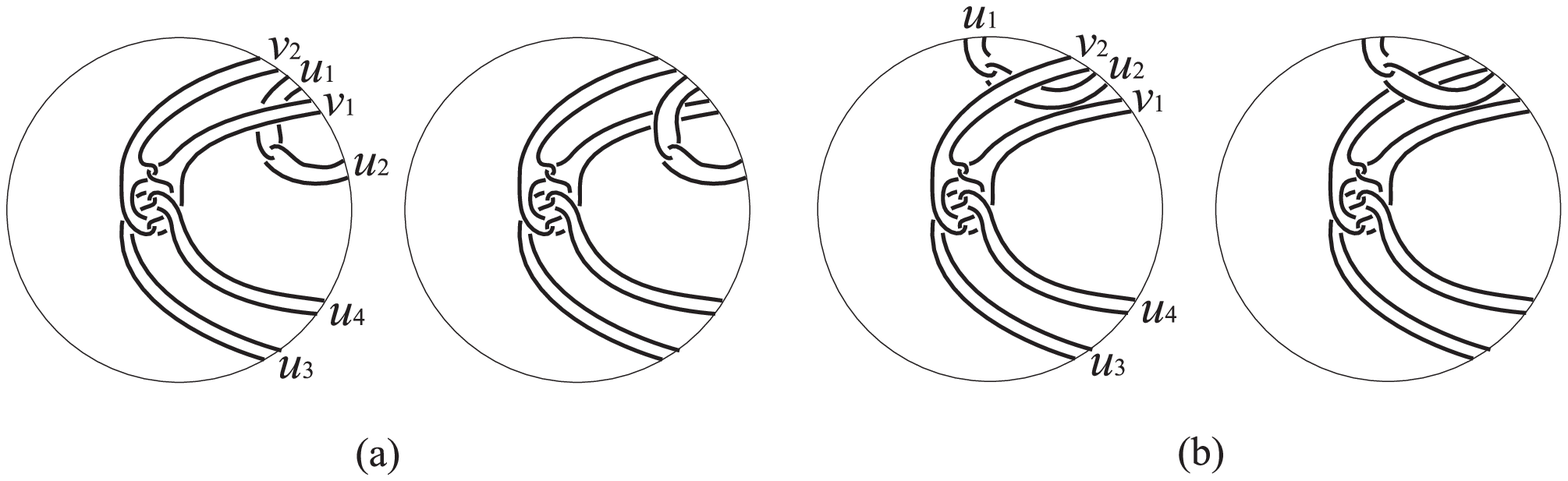}

Fig. 2.7
\end{center}

\begin{center}
\includegraphics[trim=0mm 0mm 0mm 0mm, width=.75\linewidth]
{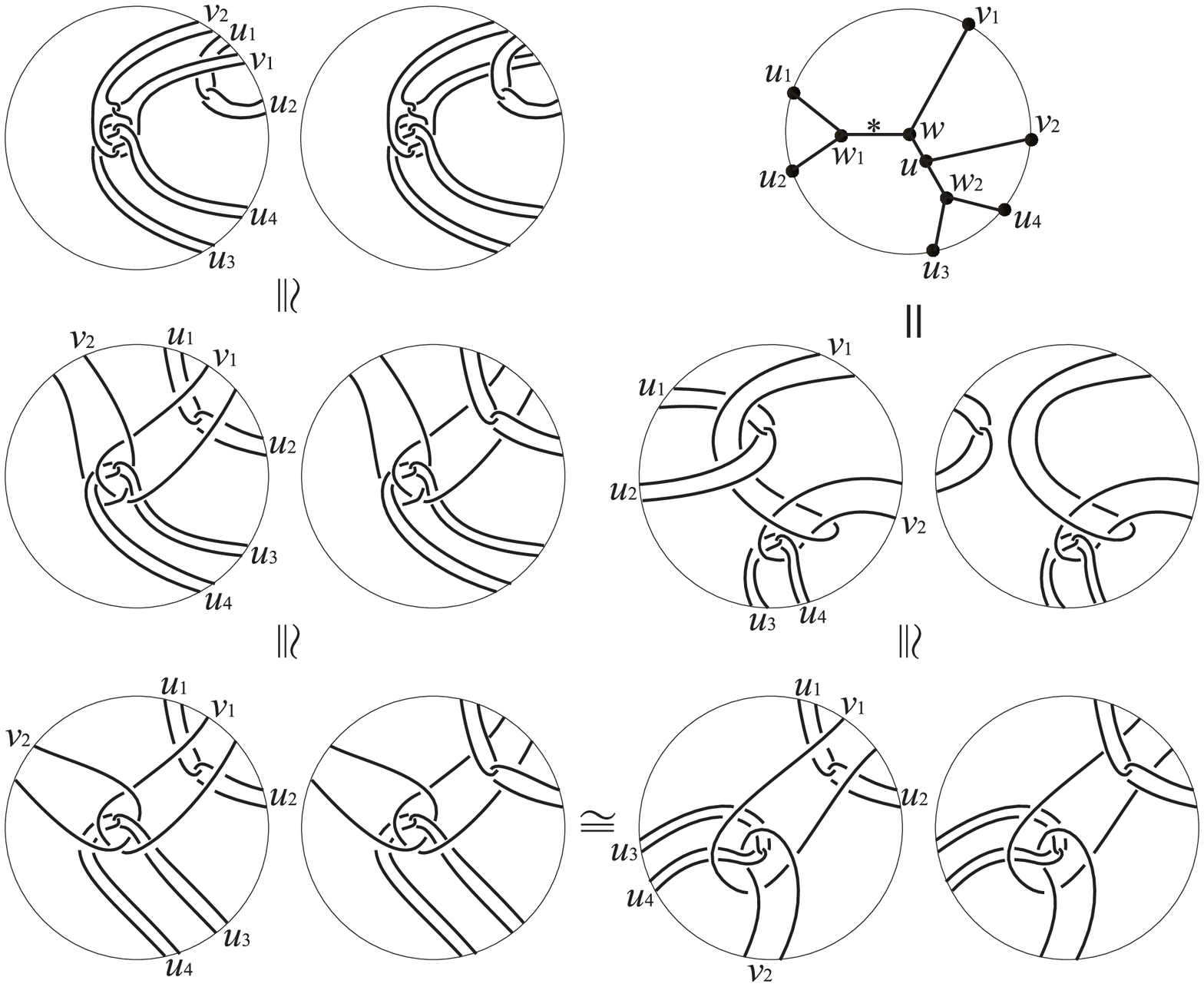}

Fig. 2.8
\end{center}

\medskip
\noindent{\bf Corollary 2.3.} (Habiro \cite{Habiro2}) 
{\it Two oriented knots $K_1$ and $K_2$ are
$C_k$-equivalent if and
only if they are related by a finite sequence of one-branched $C_k$-moves
and ambient isotopies.
$\Box$}

\medskip
Let $T_0$ be a $j$-component tangle. Let
${\cal T}(T_0)$ be the set of the tangles each element of which is homotopic
to $T_0$ relative to
$\partial B^3$. Let
${\Bbb Z}{\cal T}(T_0)$ be the free abelian group generated by ${\cal T}(T_0)$.
Let $k_1,...,k_l$ be natural numbers.
We define a {\it
singular tangle of type $(k_1,...,k_l)$} and a
subgroup ${\cal V}(k_1,...,k_l)(T_0)$ of ${\Bbb Z}{\cal T}(T_0)$ as follows.
Suppose that for each
$P\subset\{1,...,l\}$ a tangle $T_P\in{\cal T}(T_0)$ is assigned. Suppose 
that there are orientation preserving embeddings 
$h_i:B^3\rightarrow {\rm int}B^3$ $(i=1,...,l)$ such that \\
(1) $h_i(B^3)\cap h_j(B^3)\neq\emptyset$ if $i\neq j$,\\
(2) $T_P-\bigcup_{i=1}^l h_i(B^3)=T_{P'}-\bigcup_{i=1}^l h_i(B^3)$ 
for any subsets $P,P'\subset \{1,...,l\}$,\\
(3) $(h_i^{-1}(T_\emptyset),h_i^{-1}(T_{\{1,...,l\}}))$ is 
a $C_{k_i}$-move $(i=1,...,l)$, and\\
(4) $T_P\cap h_i(B^3)=\left\{
\begin{array}{ll}
T_{\{1,...,l\}}\cap h_i(B^3) & \mbox{if $i\in P$},\\
T_\emptyset\cap h_i(B^3) & \mbox{otherwise}.
\end{array}
\right.$\\
Then we call the set of tangles $\{T_P|P\subset\{1,...,l\}\}$ a {\em
singular tangle of type
$(k_1,...,k_l)$}. For a singular tangle $T=\{T_P|P\subset\{1,...,l\}\}$ of type
$(k_1,...,k_l)$, we define an element $\kappa(T)$ of ${\Bbb Z}{\cal T}(T_0)$ by 
\[\kappa(T)=\sum_{P\subset\{1,...,l\}}(-1)^{|P|}T_P.\]
Let ${\cal V}(k_1,...,k_l)(T_0)$ be the subgroup of ${\Bbb Z}{\cal T}(T_0)$
generated by all $\kappa(T)$ where $T$ varies over all
singular tangles of type $(k_1,...,k_l)$.

\medskip
\noindent{\bf Proof of Theorem 1.1.} By Corollary 2.3 it is sufficient to
show that if 
$(T_1,T_2)$ is a one-branched $C_k$-move, then $T_1-T_2\in{\cal
V}(k_1,...,k_l)(T_1)$. We first 
show that $T_1-T_2\in{\cal V}(k-k_l,k_l)(T_1)$. By Lemma 2.1 we have that
any one-branched $C_k$-move
is equivalent to a move illustrated in Fig. 2.9 (a) or (b).
Therefore we may assume that
$T_1$ and $T_2$ are as illustrated in Fig. 2.9 (a) or (b). 
Suppose that $T_1$ and $T_2$ are as in Fig. 2.9 (a) (resp. (b)). 
Let $T_3,T_4,T_5$ and $T_6$ be tangles as
illustrated in Fig. 2.10 (a) (resp. (b)). Note that $T_3$ and $T_4$ are 
ambient isotopic relative to $\partial
B^3$. Then we have that $T_1-T_2=(T_1-T_5-T_3+T_6)+(T_5-T_2-T_6+T_4)$. 
Note that $T_1$ and $T_3$,
$T_2$ and $T_4$, and $T_5$ and $T_6$ are related by a (one-branched)
$C_{k_l}$-move. Similarly $T_1$
and $T_5$, $T_3$ and $T_6$, $T_5$ and $T_2$, and $T_6$ and $T_4$ are
related by a one-branched
$C_{k-k_l}$-move. 
It is not hard to see that both $T=\{T_1,T_5,T_3,T_6\}$ and 
$T'=\{T_5,T_2,T_6,T_4\}$ are singular tangles of type 
$(k-k_l,k_l)$ and that $\varepsilon\kappa(T)=T_1-T_5-T_3+T_6$ and
$\varepsilon'\kappa(T')=T_5-T_2-T_6+T_4$, 
where $\varepsilon=\pm 1$, and $\varepsilon'=\pm 1$. 
Thus $T_1-T_2\in {\cal V}(k-k_l,k_l)(T_1)$. 

Similarly we have that if $(T'_1,T'_2)$ is a one-branched $C_{k-k_l}$-move,
then $T'_1-T'_2\in{\cal V}(k-k_l-k_{l-1},k_{l-1})(T_1)$. 
By substituting this to one-branched
$C_{k-k_l}$-moves in $B^3$ with
respect to $T_1$ and $T_5$, $T_3$ and $T_6$, $T_5$ and $T_2$, and $T_6$ and
$T_4$, we have that
$T_1-T_2\in{\cal V}(k-k_l-k_{l-1},k_{l-1},k_l)(T_1)$. Repeating this
argument we finally have the
desired conclusion. $\Box$

\begin{center}
\includegraphics[trim=0mm 0mm 0mm 0mm, width=.7\linewidth]
{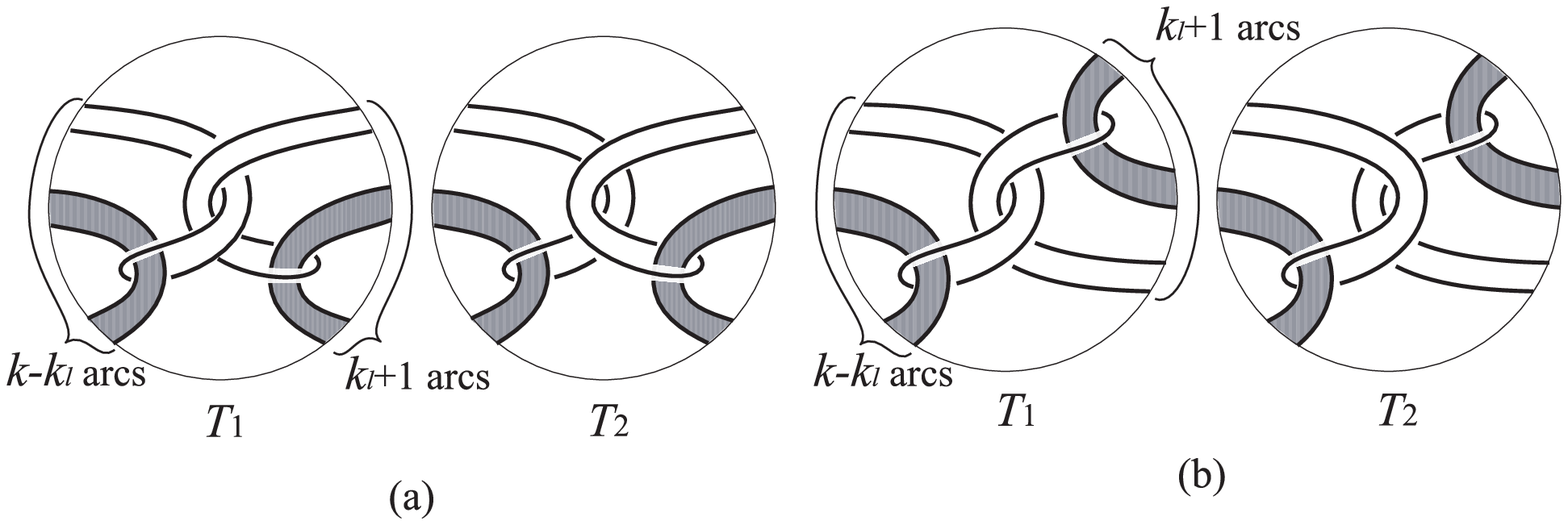}

Fig. 2.9
\end{center}

\begin{center}
\includegraphics[trim=0mm 0mm 0mm 0mm, width=.7\linewidth]
{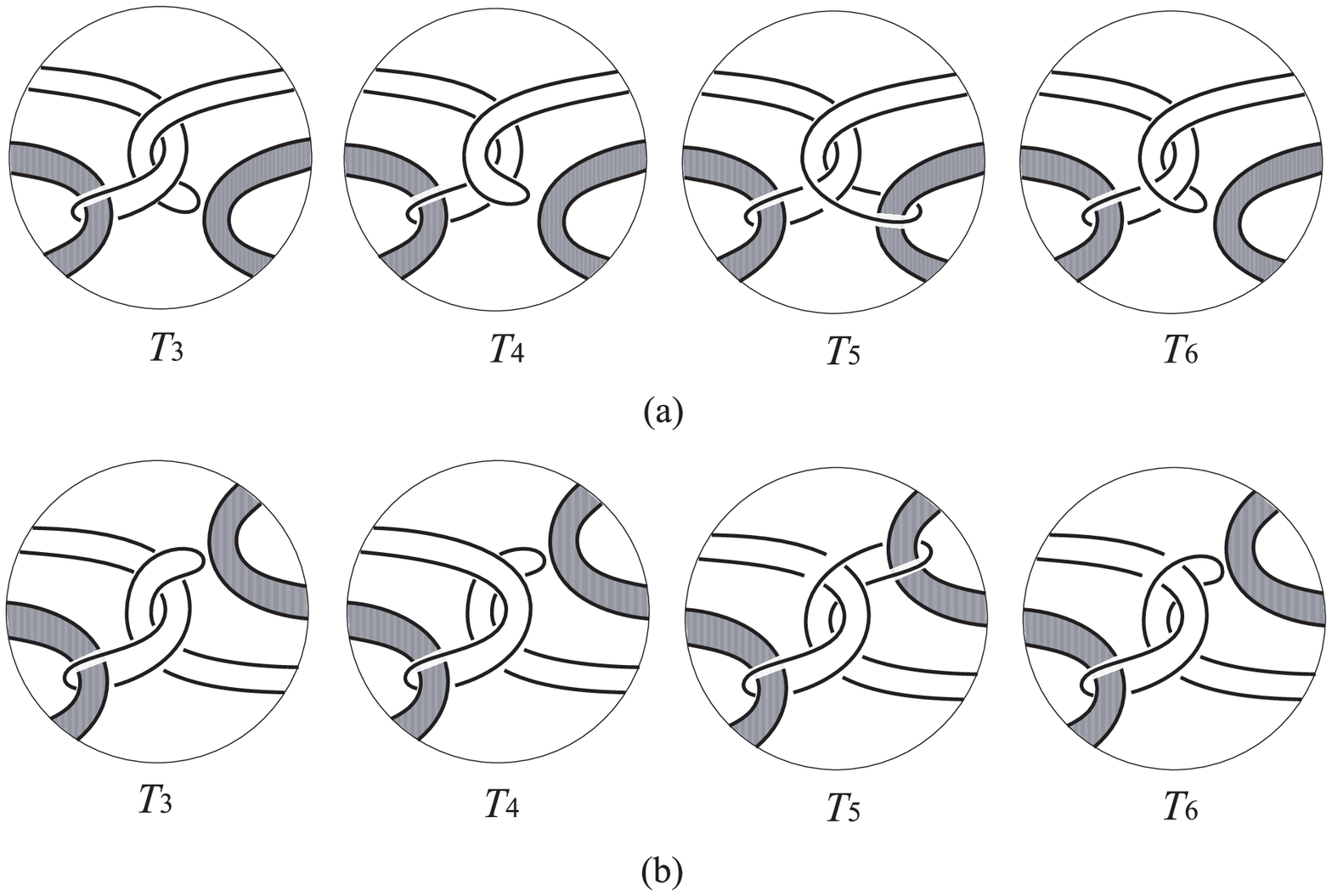}

Fig. 2.10
\end{center}

\medskip\noindent
{\bf Remark 2.4.} 
In Proof of Theorem 1.1, since $T_1$ and $T_5$, and $T_5$ and $T_2$ are 
related by a one-branched $C_{k-k_l}$-move, a one-branced $C_k$-move 
is realized by twice applications of one-branched $C_{k-k_l}$-moves. 
By Lemma 2.2, for any positive integer $k,k'$ $(k'<k)$, a $C_k$-move 
is realized by finitely many $C_{k'}$-moves. Hence $C_k$-equivalence 
implies $C_{k'}$-equivalence. 

\bigskip\noindent
{\bf 3. Band description}

A {\it $C_{1}$-link model} is a pair $(\alpha,\beta)$ where 
$\alpha$ is a disjoint union of properly embedded arcs in $B^3$ and $\beta$ 
is a disjoint union of arcs on $\partial B^3$ with 
$\partial \alpha=\partial \beta$ as illustrated in Fig. 3.1.
Suppose that a $C_{k}$-link model $(\alpha, \beta)$ is defined 
where $\alpha$ is a disjoint union of $k+1$ properly embedded arcs in $B^3$ 
and $\beta$ is a disjoint union of $k+1$ arcs on $\partial B^3$ with 
$\partial \alpha=\partial \beta$ such that $\alpha \cup \beta$ is a disjoint 
union of $k+1$ circles. Let $\gamma$ be a component of $\alpha \cup \beta$ and 
$W$ a regular neighbourhood of $\gamma$ in $(B^3-(\alpha \cup
\beta))\cup\gamma$. Let $V$ be an
oriented solid  torus, $D$ a disk in $\partial V$, $\alpha_{0}$ properly
embedded arcs in $V$ 
and $\beta_{0}$ arcs on $D$ as illustrated in Fig. 3.2. Let $\psi
:V\rightarrow W$ be an orientation
preserving homeomorphism  such that $\psi (D)=W\cap \partial B^{3}$ and 
$\psi (\alpha_{0}\cup \beta_{0})$ bounds disjoint disks in $B^3$. 
Then we call the pair $((\alpha -\gamma )\cup \psi (\alpha_{0}), 
(\beta -\gamma )\cup\psi (\beta_{0}))$ a {\it $C_{k+1}$-link model.} 
We also say that the pair $((\alpha -\gamma )\cup \psi (\alpha_{0}), 
(\beta -\gamma )\cup\psi (\beta_{0}))$ is a {\it double} of
$(\alpha,\beta)$ with respect to the
component $\gamma$.  A {\it link model} is a 
$C_{k}$-link model for some $k$. 

\begin{center} 
\includegraphics[trim=0mm 0mm 0mm 0mm, width=.17\linewidth]
{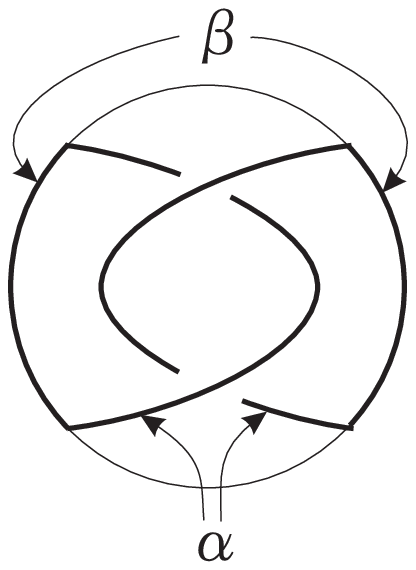}

Fig. 3.1
\end{center}

\begin{center}
\includegraphics[trim=0mm 0mm 0mm 0mm, width=.2\linewidth]
{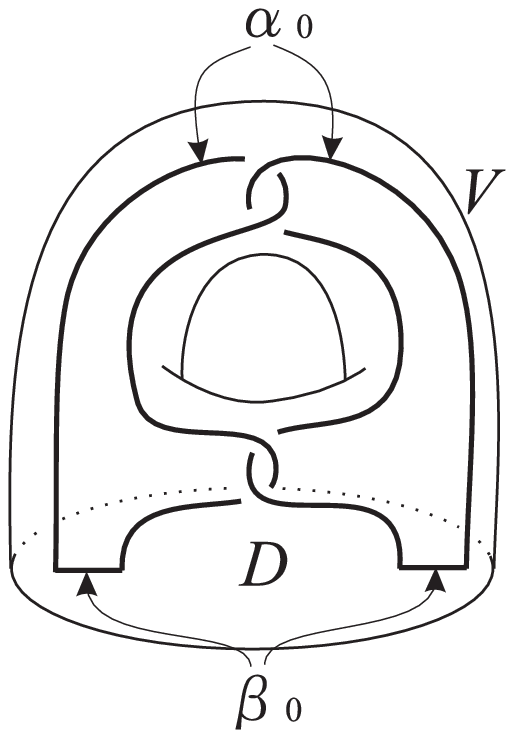}

Fig. 3.2
\end{center}

Let $(\alpha_{1},\beta_{1}),...,(\alpha_{l},\beta_{l})$ be link 
models. Let $K$ be an oriented knot (resp. a tangle). Let
$\psi_{i}:B^{3}\rightarrow
S^3$ (resp. $\psi_{i}:B^{3}\rightarrow
{\rm int}B^3$)  be an orientation preserving embedding for $i=1,...,
l$ and 
$b_{1,1},b_{1,2},...,b_{1,\rho(1)},b_{2,1},b_{2,2},...,b_{2,\rho(2)},
...,b_{l,1},b_{l,2},...,b_{l,\rho(l)}$
mutually disjoint disks embedded in $S^3$ (resp. $B^3$). Suppose that 
they satisfy the following
conditions;\\
(1) $\psi_{i}(B^{3})\cap \psi_{j}(B^{3})=\emptyset$ if $i\neq j$,\\
(2) $\psi_{i}(B^{3})\cap K=\emptyset$  for each $i$,\\
(3) $b_{i,k}\cap K=\partial b_{i,k}\cap K$ is an arc for each $i,k$,\\
(4) $b_{i,k}\cap (\bigcup_{j=1}^{l} \psi_{j}(B^{3}))=
\partial b_{i,k}\cap \psi_{i}(B^{3})$ is a component of 
$\psi_{i}(\beta_{i})$ for each $i,k$,\\
(5) ($\bigcup_{k=1}^{\rho(i)}b_{i,k})\cap \psi_{i}(B^{3})
=\psi_{i}(\beta_{i})$ for each $i$.\\
Let $J$ be an oriented knot (resp. a tangle) defined by
\[
J=K\cup (\bigcup_{i,k}\partial b_{i,k})\cup 
(\bigcup_{i=1}^{l}\psi_{i}(\alpha_{i})) - 
\bigcup_{i,k}{\rm int}(\partial b_{i,k}\cap K) - 
\bigcup_{i=1}^{l}\psi_{i}({\rm int}\beta_{i}),
\]
where the orientation of $J$ 
coincides that of $K$ on $K-\bigcup_{i,k}b_{i,k}$ 
if $K$ is oriented. Then we say that 
$J$ is a {\it band sum} of $K$ and link models
$(\alpha_{1},\beta_{1}),...,
(\alpha_{l},\beta_{l})$. We call each $b_{i,k}$ a {\it band}. 
Each image $\psi_{i}(B^{3})$ is called a {\it link ball}. In particular if
$(\alpha_i,\beta_i)$ is a $C_k$-link
model then $b_{i,k}$ is called a {\it $C_k$-band} and $\psi_{i}(B^{3})$ is
called a {\it $C_k$-link
ball}. We set
${\cal B}_i=((\alpha_i,\beta_i),\psi_i,\{b_{i,1},...,b_{i,\rho(i)}\})$
and call ${\cal B}_i$ a {\it
chord}. In particular ${\cal B}_i$ is called a {\it $C_k$-chord} when
$(\alpha_i,\beta_i)$ is a
$C_k$-link model. We denote $J$  by
$J=\Omega(K;\{{\cal B}_1,...,{\cal B}_l\})$ and call it a {\it band
description} of $J$. We also say that
$J$ is a band sum of
$K$ and chords ${\cal B}_1,...,{\cal B}_l$.

\medskip
\noindent{\bf Sublemma 3.1.} {\it Let $(T_1,T_2)$ be a $C_k$-move. Then
there is a $C_k$-link model
$(\alpha,\beta)$ such that $(T_1,T_2)\cong(\alpha,\hat{\beta})$ where
$\hat{\beta}$ is a
slight push in of $\beta$.}

\medskip
\noindent{\bf Proof.} It is clearly true for $k=1$. Suppose
that $(T_1,T_2)$ is a double of a
$C_{k-1}$-move $(U_1,U_2)$ with respect to the components $u_1$ and $u_2$
and $(\alpha',\beta')$ is a
$C_{k-1}$-link model such that
$(U_1,U_2)\cong(\alpha',\hat{\beta'})$. 
Then by the deformation illustrated in Fig. 3.3, we have that 
there is a double $(\alpha,\beta)$ of $(\alpha',\beta')$ with 
respect to the component which corresponds to $u_1$ and $u_2$ 
such that $(\alpha,\hat{\beta})$ 
is equivalent to a double of $(\alpha',\hat{\beta'})$ and therefore
equivalent to $(T_1,T_2)$. $\Box$

\begin{center}
\includegraphics[trim=0mm 0mm 0mm 0mm, width=.7\linewidth]
{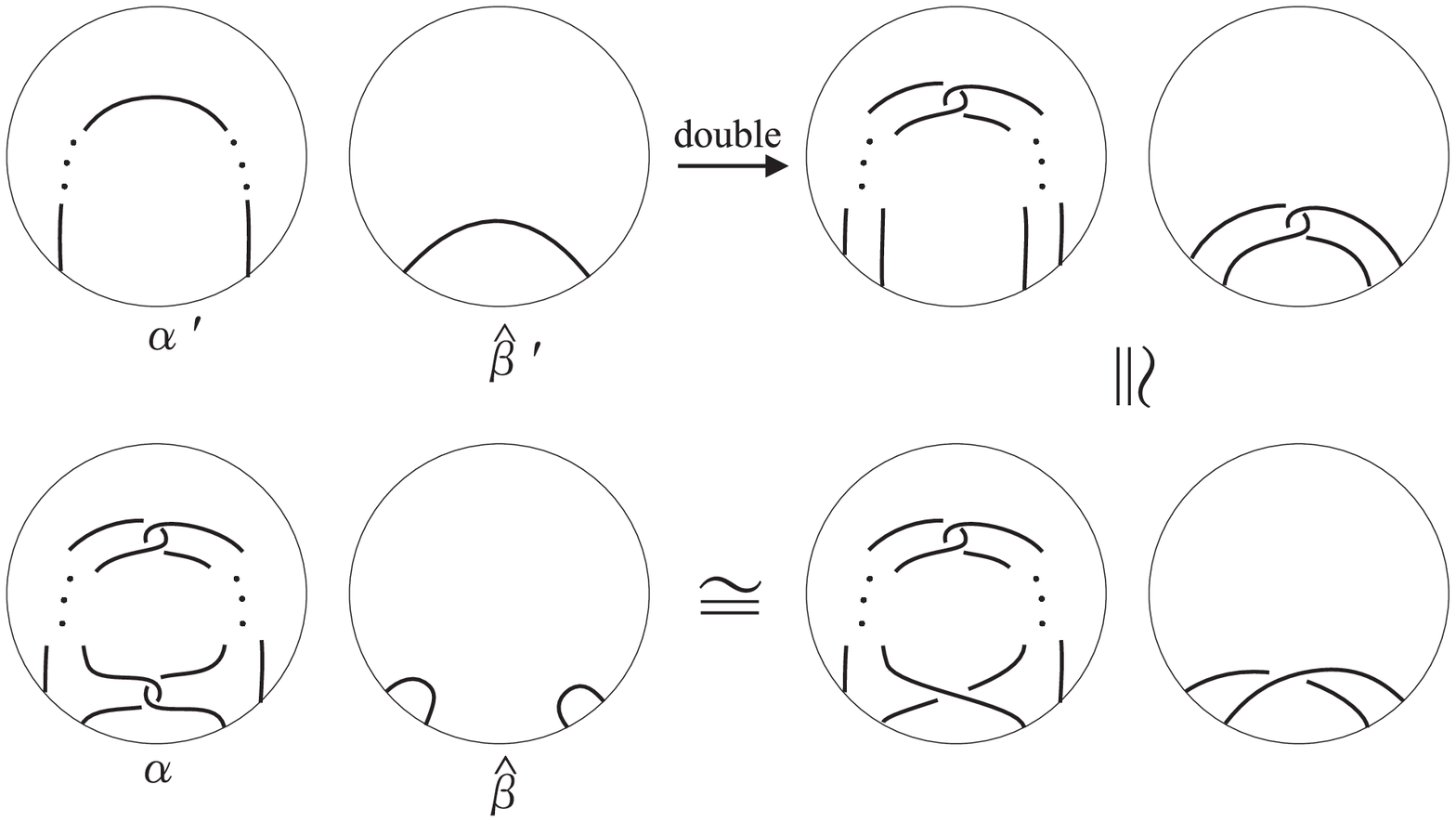}

Fig. 3.3
\end{center}

Fig. 3.3 indicates that for a double $(\alpha,\beta)$ of 
$(\alpha',\beta')$, $(\alpha,\hat{\beta})$ is equivalent 
to a double of $(\alpha',\hat{\beta'})$. 
Thus we have 

\medskip\noindent
{\bf Sublemma 3.2.} {\it If $(\alpha,\beta)$ is a $C_k$-link 
model, then $(\alpha,\hat{\beta})$ is equivalent to a $C_k$-move. 
$\Box$}

\medskip
\noindent{\bf Sublemma 3.3.} {\it Let $(T_1,T_2)$ be a $C_k$-move.
Then there is a
$C_k$-link model
$(\alpha,\beta)$ such that a band sum of $T_1$
and $(\alpha,\beta)$ is ambient isotopic to 
$T_2$ relative to $\partial B^3$.}

\medskip
\noindent{\bf Proof.} Since the inverse $(T_2,T_1)$ is also a
$C_k$-move we have by Sublemma 3.1 that there is a
$C_k$-link model
$(\alpha,\beta)$ such that $(T_2,T_1)\cong(\alpha,\hat{\beta})$. It is easy
to see that $\alpha$ is a band sum of $\hat{\beta}$ and
$(\alpha,\beta)$ up to ambient isotopy relative to $\partial B^3$. 
Therefore we have the result. $\Box$

\medskip From now on we consider knots up to ambient isotopy of 
$S^3$ and tangles up to ambient isotopy of
$B^3$ relative to $\partial B^3$ without explicit mention. As an immeditate
consequence of Sublemmas 3.2 and 3.3 we have 

\medskip
\noindent{\bf Sublemma 3.4.} {\it Let $K$ and $J$ be oriented
knots. Then $K$ and $J$ are related by a $C_k$-move if and 
only if $J$ is a band sum of 
$K$ and a $C_k$-link model. $\Box$}

\medskip
Let $K$ be a knot and $J=\Omega(K;\{{\cal B}_1,...,{\cal B}_{l}\})$ 
a band sum of $K$ and $C_{k_i}$-chords ${\cal B}_i$ $(i=1,...,l)$. 
We define an element $\kappa(J)$ of ${\Bbb Z}{\cal K}$ by 
\[\kappa(J)=\sum_{P\subset\{1,...,l\}}(-1)^{|
P|} \Omega\left(K;\bigcup_{i\in P}\{{\cal B}_{i}\}\right).\]
By Subemmas 3.2 and 3.3 we have that the subgroup 
${\cal V}(k_1,...,k_l)$ of ${\Bbb Z}{\cal K}$ is 
generated by all $\kappa(J)$ where $J$ varies over all
band sums of knots and $C_{k_i}$-chords $(i=1,...,l)$.

\medskip
\noindent{\bf Sublemma 3.5.} {\it Let $K$, $J$ and $I$ be oriented
knots. Suppose that
$J=\Omega(K;\{{\cal B}_1,...,{\cal B}_l\})$ for some chords ${\cal
B}_1,...,{\cal B}_l$
and $I=\Omega(J;\{{\cal B}\})$ for some $C_k$-chord ${\cal B}$. Then there
is a $C_k$-chord ${\cal B}'$ such that 
$I=\Omega(K;\{{\cal B}_1,...,{\cal B}_l,{\cal B}'\})$. 
Moreover, if there is a subset $P$ of $\{1,...,l\}$ such that 
the link ball and the bands of $\cal B$ intersect neither the link ball 
nor the bands of ${\cal B}_j$ for any $j\in\{1,...,l\}\setminus P$, 
then $\Omega(\Omega(K;\bigcup_{i\in P}\{{\cal B}_i\});\{{\cal B}\})=
\Omega(K;(\bigcup_{i\in P}\{{\cal B}_i\})\cup\{{\cal B}'\})$.}

\medskip
\noindent{\bf Proof.} If the bands and the
link ball of ${\cal B}$ are disjoint from those of ${\cal B}_1,...,{\cal
B}_l$ then we have
that
$I=\Omega(K;\{{\cal B}_1,...,{\cal B}_l,{\cal B}\})$. If not then we
deform $I$ up to
ambient isotopy as follows. First we thin and shrink the bands and the link
ball of ${\cal B}$
so that they are thin enough and small enough respectively. If the link
ball of ${\cal B}$
intersects the bands and the link balls of ${\cal B}_1,...,{\cal
B}_l$ then we move the link
ball of ${\cal B}$ so that they does not intersect. Then we slide the
bands of
${\cal B}$ along $J$ so that the intersection of the bands with $J$ is
disjoint from the
bands and the link balls of ${\cal B}_1,...,{\cal B}_l$. Then we
sweep the bands of
${\cal B}$ out of the link balls of ${\cal B}_1,...,{\cal B}_l$. Note
that this is
always possible since the tangles are trivial. Finally we sweep the
intersection of the bands of
${\cal B}$ and the bands of ${\cal B}_1,...,{\cal B}_l$ out of the
intersection of the bands
of ${\cal B}_1,...,{\cal B}_l$ and $K$. Let ${\cal B}'$ be the result
of the deformation of
${\cal B}$ described above. Then it is not hard to see that 
${\cal B}'$ is a desired chord. $\Box$

\medskip
By repeated applications of Sublemmas 3.4 and 3.5 we immediately 
have the following lemma.

\medskip
\noindent{\bf Lemma 3.6.} {\it Let $k$ be a positive integer 
and let $K$ and $J$ be oriented knots. 
Then $K$ and $J$ are $C_{k}$-equivalent if and only if $J$ is a band
sum of $K$ and some $C_{k}$-link models. $\Box$}

\medskip
As in the definition of $C_k$-move we define iteratedly doubled strings as
follows. A {\it $0$-double pattern} is $\{{\bf o}\}\times[0,1]$ in 
$B^2\times [0,1]$ where ${\bf o}$ is the center of
$B^2$. Suppose that a $k$-double pattern $A$ in $B^2\times [0,1]$ is
defined. Let $N$ be a regular
neighbourhood of a component $\gamma$ of $A$ in $(B^2\times
[0,1]-A)\cup\gamma$ such that
$N\cap(\partial B^2\times[0,1])=\emptyset$. Let
$\psi:B^2\times[0,1] \rightarrow N$ be a homeomorphism with
$\psi(B^2\times\{0,1\})=N\cap(B^2\times\{0,1\})$. Let
$\alpha$ be a disjoint union of properly embedded arcs in $B^2\times[0,1]$
as illustrated in Fig. 1.1. Then $(A-\gamma)\cup\psi(\alpha)$ 
is called a {\it $(k+1)$-double pattern}. Let $N$ be a regular
neighbourhood of a properly embedded arc in $B^3$. Let
$\psi:B^2\times[0,1]\rightarrow N$ be a
homeomorphism. Then the image $\psi(A)$ of a $k$-double pattern $A$ is
called a {\it $k$-double
string}. Note that a \lq crossing change' between $k$-double 
strng and $j$-double string is equivalent to a $C_{k+j+1}$-move. 

\medskip
\noindent{\bf Sublemma 3.7.} {\it Let $(\alpha,\beta)$ be a
$C_k$-link model and $\gamma$ a
component of $\alpha\cup\beta$. Let $D$ be a disk in $B^3$ such that
$\partial D=\gamma$ and
${\rm int}D\cap \partial B^3=\emptyset$. Let
$\delta$ be a properly embedded unknotted arc in $B^3$ that intersects $D$
transversally at one point
in ${\rm int}D$. Let $N$ be a regular neighbourhood of $\delta$ in
$B^3-\gamma$. Then there exist a
$(k-1)$-double string 
${\cal D}$ in $N$ and an orientation preserving homeomorphism
$\varphi:B^3\rightarrow B^3$ such that 
$\varphi|_{\gamma\cap\beta}={\rm id}|_{\gamma\cap\beta}$, 
and $\varphi(\alpha\cup\gamma)={\cal D}\cup\gamma$. }

\medskip\noindent
In the lemma above, we note that ${\cal D}\cup(\alpha\cap \gamma)$ 
is a $k$-double string. 

\medskip
\noindent{\bf Proof.} The case $k=1$ is clear. Suppose that it
is shown for $k-1$.
Let $(\alpha,\beta)$ be a double of a $C_{k-1}$-link model
$(\alpha',\beta')$ with respect
to the component $\gamma'$ of $\alpha'\cup\beta'$.

If $\gamma$ is already a component of
$\alpha'\cup\beta'$, then by the deformation of $\alpha$ in a regular
neighbourhood of $\gamma'$ as
illustrated in Fig. 3.4 and by the assumption, we have the result.

Next suppose that $\gamma$ is not a component of $\alpha'\cup\beta'$. In
other words $\gamma$ is
contained in a regular neighbourhood $W$ of $\gamma'$. Let $D'$ be a disk
in $B^3$ such that $\partial
D'=\gamma'$ and
${\rm int}D'\cap \partial B^3=\emptyset$. Let
$\delta'$ be a properly embedded unknotted arc in $B^3$ that intersects
$D'$ transversally at one
point in ${\rm int}D'$. Let $N'$ be a regular neighbourhood of $\delta'$ in
$B^3-\gamma'$. Then
by the assuption there is a
$(k-2)$-double string ${\cal D}'$ in
$N'$ and an
orientation preserving homeomorphism $\varphi':B^3\rightarrow B^3$ such that
$\varphi'|_{\gamma'\cap\beta'}={\rm id}|_{\gamma'\cap\beta'}$ and
$\varphi'(\alpha'\cup\gamma')={\cal D}'\cup\gamma'$. By modifying $\varphi'$
if necessary we may suppose that
$\varphi'|_{W\cap\partial B^3}={\rm id}|_{W\cap\partial B^3}$. Then
we have that
$\varphi'(\alpha\cup\gamma)$ is as illustrated in Fig. 3.5 (a). Then by the
deformation illustrated
in Fig. 3.5 we have the result. $\Box$

\begin{center}
\includegraphics[trim=0mm 0mm 0mm 0mm, width=.4\linewidth]
{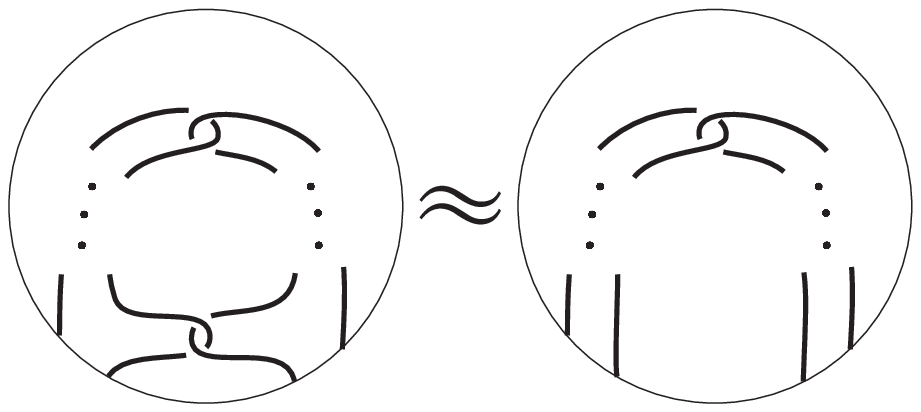}

Fig. 3.4
\end{center}

\begin{center}
\includegraphics[trim=0mm 0mm 0mm 0mm, width=.7\linewidth]
{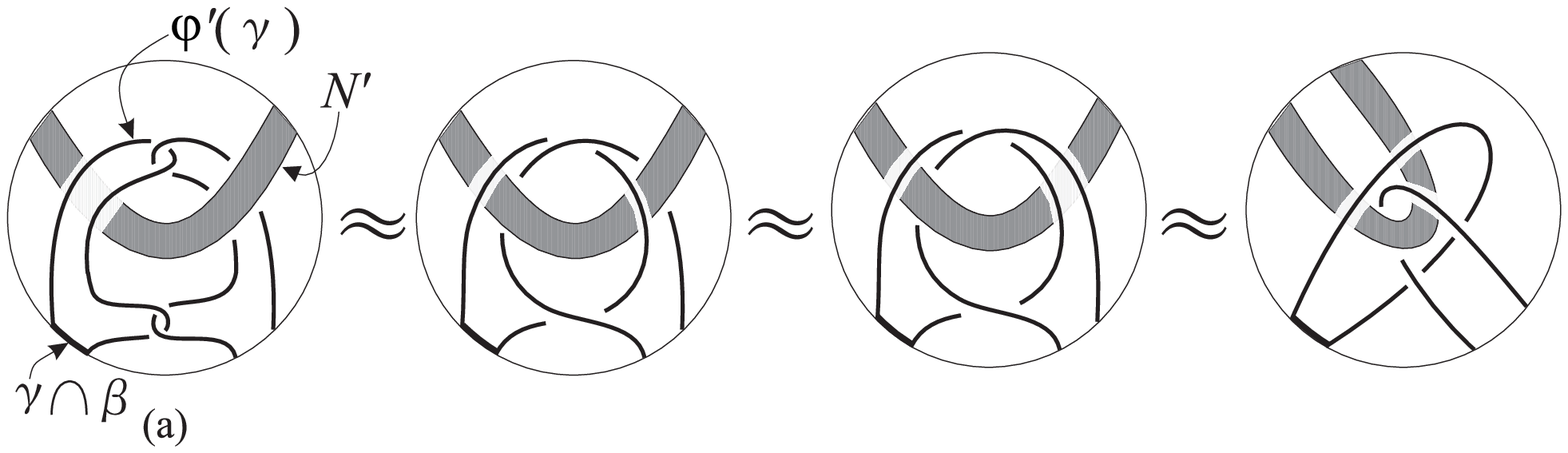}

Fig. 3.5
\end{center}

\medskip
\noindent{\bf Lemma 3.8.} {\it 
 Let $K$, $J=\Omega(K;\{{\cal B}_1,...,{\cal B}_l,{\cal B}_0\})$ 
 and $I=\Omega(K;\{{\cal B}_1,...,{\cal B}_l,{\cal B}'_0\})$ be 
 oriented knots, where 
${\cal B}_1,...,{\cal B}_l$ are chords and ${\cal B}_0,{\cal B}'_0$ 
are $C_k$-chords. Suppose that $J$ and $I$ differ locally as illustrated 
in Fig $3.6$ {\rm (a), (b)}, i.e., 
$I$ is obtained from $J$ by a crossing change 
between $K$ and a band of ${\cal B}_0$. Then $J$ and $I$ are related 
by a $C_{k+1}$-move. Moreover, there is a $C_{k+1}$-chord $\cal B$ such 
that $\Omega(K;(\bigcup_{i\in P}\{{\cal B}_i\})\cup\{{\cal B}_0\})=
\Omega(K;(\bigcup_{i\in P}\{{\cal B}_i\})\cup\{{\cal B}'_0,{\cal B}\})$ 
for any subset $P$ of $\{1,...,l\}$. }

\medskip
\noindent{\bf Proof.} By shrinking the band and pulling the
link ball as illustrated in 
Fig. 3.6 it is sufficient to show the case that $K$ is near the link
ball. Then by Sublemma 3.7 
we deform the strings in the link ball without disturbing a neighbourhood
of the band. Then the crossing change is realized by a 
$C_{k+1}$-move. See Fig. 3.7. In Fig. 3.7, there is a $3$-ball $B$ in 
$S^3$ and a homeomorphism $h:B\rightarrow B^3$ 
such that $(h(J), h(I))$ is a $C_{k+1}$-move and $B$ is 
disjoint from the link ball and the bands of any chord ${\cal B}_i$ 
$(i=1,...,l)$. Thus by Sublemmas 3.3 and 3.5, we have the 
latter assertion. $\Box$

\begin{center}
\includegraphics[trim=0mm 0mm 0mm 0mm, width=.7\linewidth]
{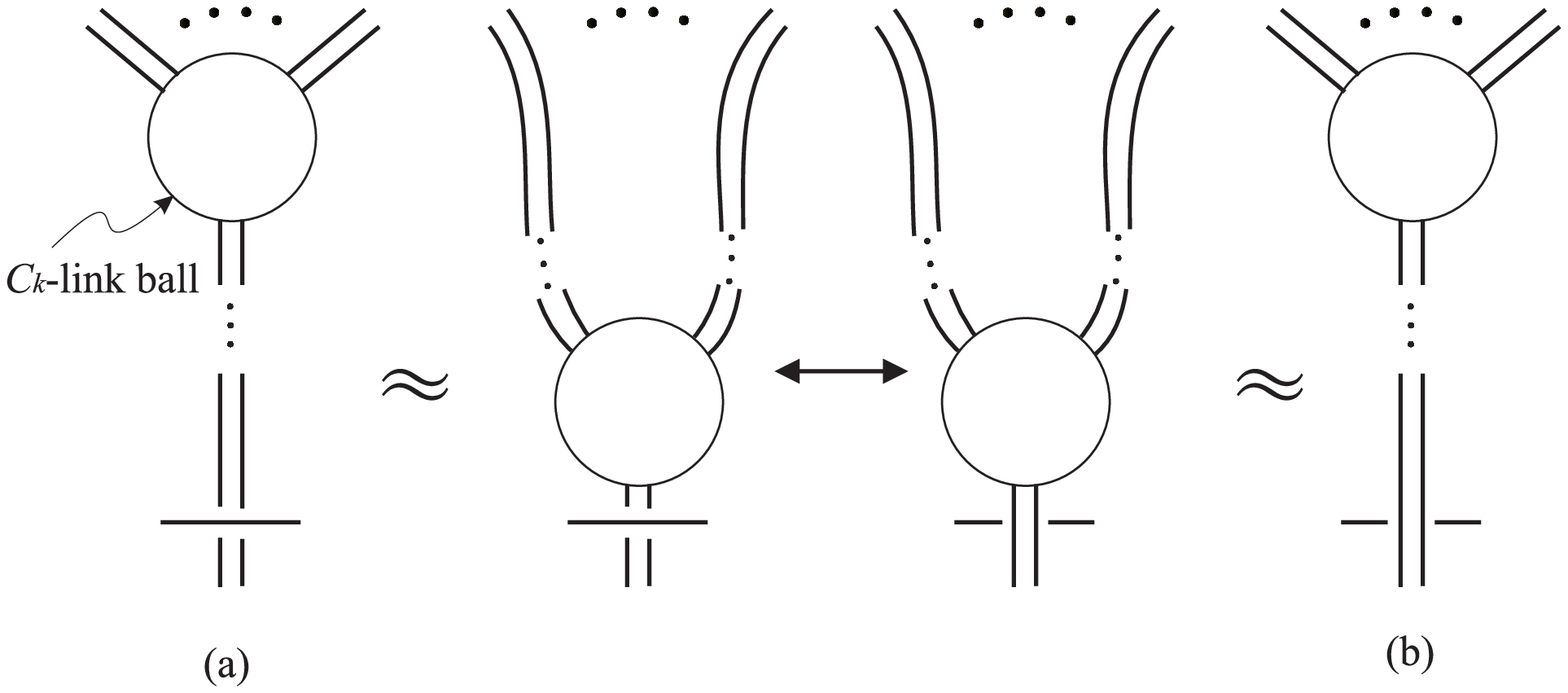}

Fig. 3.6
\end{center}

\begin{center}
\includegraphics[trim=0mm 0mm 0mm 0mm, width=.4\linewidth]
{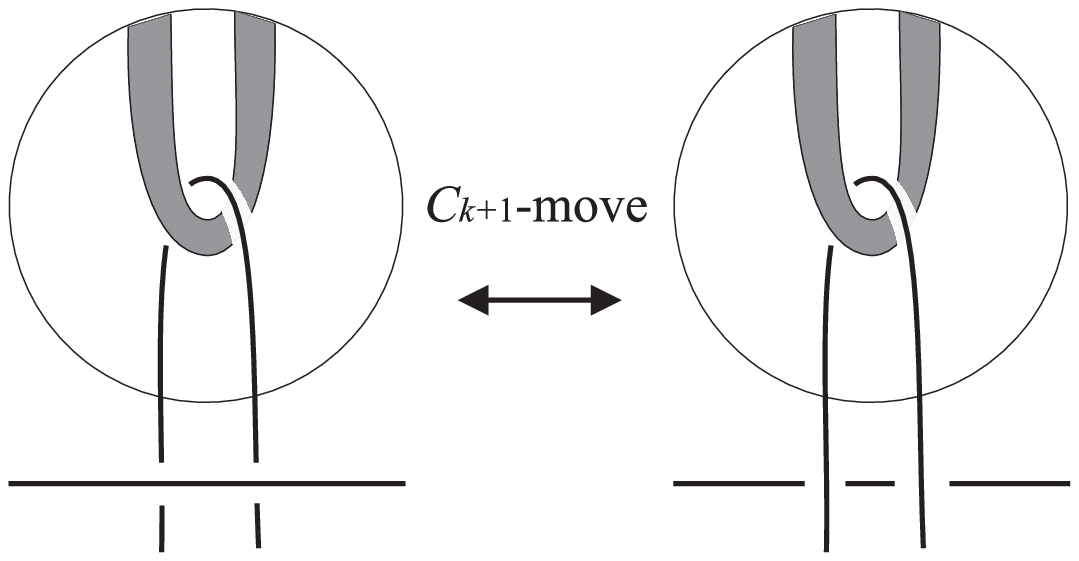}

Fig. 3.7
\end{center}

\medskip
\noindent{\bf Lemma 3.9.} {\it  Let $K$, 
$J=\Omega(K;\{{\cal B}_1,...,{\cal B}_l,{\cal B}_{0j}, 
{\cal B}_{0k}\})$ and $I=\Omega(K;\{{\cal B}_1,...,{\cal B}_l,
{\cal B}'_{0j},{\cal B}'_{0k}\})$ be oriented knots, where 
${\cal B}_1,...,{\cal B}_l$ are chords and 
${\cal B}_{0j},{\cal B}'_{0j}$ $($resp. ${\cal B}_{0k},{\cal B}'_{0k})$ 
are $C_j$-chords $($resp. $C_k$-chords$)$. 
Suppose that $J$ and $I$ differ locally as illustrated in Fig. {\rm3.8}. 
Then $J$ and $I$ are related by a $C_{j+k}$-move. 
Moreover, there is a $C_{j+k}$-chord $\cal B$ such 
that $\Omega(K;(\bigcup_{i\in P}\{{\cal B}_i\})\cup\{{\cal B}_{0j}, 
{\cal B}_{0k}\})=\Omega(K;(\bigcup_{i\in P}\{{\cal B}_i\})\cup
\{{\cal B}'_{0j},{\cal B}'_{0k},{\cal B}\})$ 
for any subset $P$ of $\{1,...,l\}$. }

\medskip\noindent
We call the change from $J$ to $I$ in Lemma 3.9 a {\it band exchange}.

\begin{center}
\includegraphics[trim=0mm 0mm 0mm 0mm, width=.7\linewidth]
{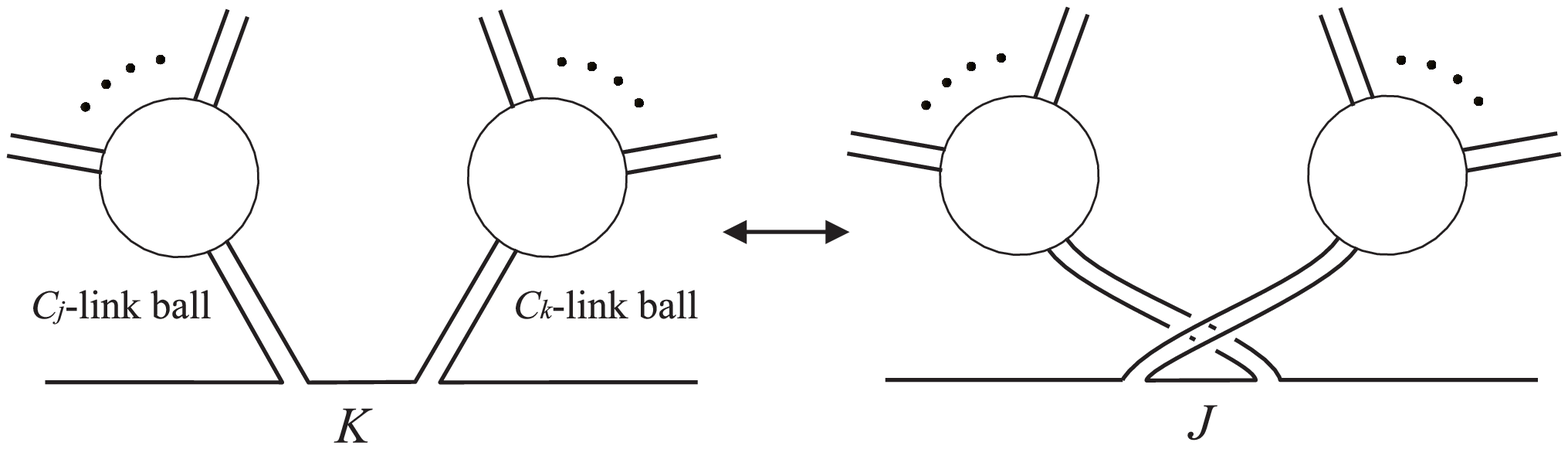}

Fig. 3.8
\end{center}

\medskip
\noindent{\bf Proof.} First we deform the strings in the link
balls as stated in Sublemma 3.7 then we slide one 
of the two bands along the 
other band and then perform a $C_{j+k}$-move. See Fig. 3.9. 
In Fig. 3.9, there is a $3$-ball $B$ in 
$S^3$ and a homeomorphism $h:B\rightarrow B^3$ 
such that $(h(J), h(I))$ is a $C_{j+k}$-move and $B$ is 
disjoint from the link ball and the bands of any chord ${\cal B}_i$ 
$(i=1,...,l)$. Thus by Sublemmas 3.3 and 3.5, we have the 
latter assertion. $\Box$

\begin{center}
\includegraphics[trim=0mm 0mm 0mm 0mm, width=.7\linewidth]
{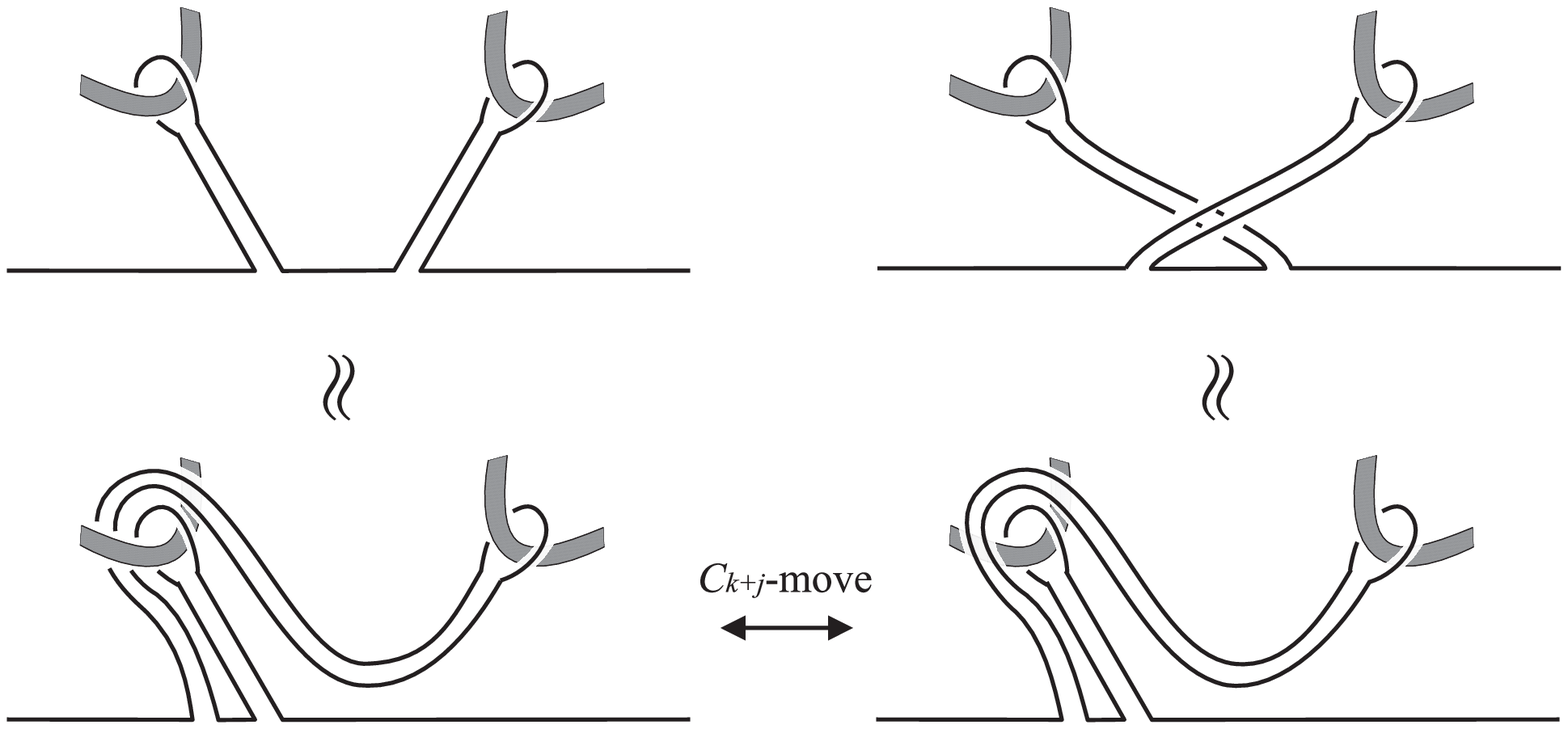}

Fig. 3.9
\end{center}

\medskip
\noindent
{\bf Proof of Theorem 1.2.} 
Let $k_1,...,k_l$ $(l\geq 2)$ be positive integer and 
$k=k_1+\cdots+ k_l$. 
Let $K_0$ be a knot and $K_1$ a band sum of $K_0$ and $C_{k_j}$-chords 
${\cal B}_{k_j,j}$ $(j=1,...,l)$. 
It is sufficient to show that
\[\mu_k\left(\sum_{P\subset\{1,...,l\}}(-1)^{|
P|} \Omega\left(K_0;\bigcup_{j\in P}\{{\cal B}_{k_j,j}\}\right)\right)
=0\in {\Bbb Z}{\cal K}/({\cal V}(k)+{\cal R}_\#),\]
where $\mu_k:{\Bbb Z}{\cal K}\rightarrow
{\Bbb Z}{\cal K}/({\cal V}(k)+{\cal R}_\#)$ is the quotient 
homomorphism. 

Set
\[K_P=\Omega\left(K_0;\bigcup_{j\in P}\{{\cal
B}_{k_j,j}\}\right).\]
{\em In the following we consider knots up to
$C_{k}$-equivalence and by a symbol 
$K_1$ $($resp. $K_P)$ we express a knot
that is $C_{k}$-equivalent to
$K_1$ $($resp. $K_P)$.} We will deform the form of band
description of $K_1$ step by step. At each step
$K_1$ is expressed as a band sum of $K_0$ and some chords such that 
each $K_P$ is a band sum of
$K_0$ and some subset of the chords of $K_1$. To be more precise
\[K_1=\Omega\left(K_0;\bigcup_{i,j}\{{\cal B}_{i,j}\}\right)\]
at each step where ${\cal B}_{i,j}$ is a
$C_i$-chord for some $i$ with $1\leq i<k$ and it has an 
associated subset $\omega({\cal B}_{i,j})\subset\{1,...,l\}$ with 
$\sum_{t\in\omega({\cal B}_{i,j})}k_t\leq i$ such that for each 
$P\subset\{1,...,l\}$
\[(*)\ \ K_P=\Omega\left(K_0;\bigcup_{\omega({\cal B}_{i,j})
\subset P}\{{\cal B}_{i,j}\}\right).\]
A chord ${\cal B}_{i,j}$ is called a {\it local chord} if there is a 3-ball
$B$ such that $B$ contains all of the bands and the link ball of 
${\cal B}_{i,j}$, $B$ does not intersect any other bands and link
balls, and that $(B,B\cap K_0)$ is a trivial ball-arc pair. Such a local
chord ${\cal B}_{ij}$ represents a knot $K_{ij}$ connected summed to $K_0$. 
The final goal of the step by step deformation is a band sum of $K_0$
and some local chords ${\cal B}_{ij}$'s so that 
$K_1$ is a connected sum of $K_0$ and $K_{ij}$'s. 
Since $\mu_k(K\#K')=\mu_k(K+K')\in 
{\Bbb Z}{\cal K}/({\cal V}(k)+{\cal R}_\#)$, 
we have 
\[\begin{array}{rl}
 &\displaystyle{\mu_k\left(
 \sum_{P\subset\{1,...,l\}}
 (-1)^{| P|}\Omega\left(K_0;\bigcup_{\omega({\cal B}_{i,j})\subset P}
 \{{\cal B}_{i,j}\}\right)\right)}\\
 =&\displaystyle{\mu_k\left(
 \sum_{P\subset\{1,...,l\}}(-1)^{| P|}\left(K_0
 +\sum_{\omega({\cal B}_{i,j})\subset P}K_{i,j}\right)\right)}\\
=&\displaystyle{\mu_k\left(
 \sum_{P\subset\{1,...,l\}}(-1)^{| P|}K_0+
 \sum_{P\subset\{1,...,l\}}(-1)^{| P|}\left(\sum_{\omega
 ({\cal B}_{i,j})\subset P}K_{i,j}\right)\right)}\\
=&\displaystyle{\mu_k\left(
 0+\sum_{i,j}\left(\sum_{P\subset\{1,...,l\},\omega({\cal B}_{i,j})
 \subset P}(-1)^{| P|}K_{i,j}\right)\right)}.
\end{array}\]
We consider the coefficient of $K_{i,j}$. 
Since $\sum_{t\in\omega({\cal B}_{i,j})}k_t< k$, 
$\omega({\cal B}_{i,j})$ is a proper subset of
$\{1,...,l\}$. We may assume that $\omega({\cal B}_{i,j})$ 
does not contain $a\in\{1,...,l\}$.
Then we have that
\[\begin{array}{rl}
\displaystyle{ 
\sum_{P\subset\{1,...,l\},\omega({\cal B}_{i,j})\subset P}
(-1)^{| P|}}=
&\displaystyle{
\sum_{P\subset\{1,...,l\}\setminus\{a\},\omega({\cal B}_{i,j})
\subset P}(-1)^{| P|}}\\
&\hspace*{3em}+
\displaystyle{\sum_{P\subset\{1,...,l\}\setminus\{a\},\omega({\cal
B}_{i,j})\subset P}(-1)^{| P\cup\{a\}|}=0.}
\end{array}\]
Thus, we have the conclusion if we can get a desired band 
description.  

Now we will deform the band description of $K_1$ into a desired form. 
We first set $\omega({\cal B}_{k_j,j})=\{j\}$ for $j=1,...,l$. 
Then we have $\sum_{t\in\omega({\cal B}_{k_j,j})}k_t=k_j$ and 
\[K_P=\Omega\left(K_0;\bigcup_{\omega({\cal
B}_{k_j,j})\subset P}\{{\cal B}_{k_j,j}\}\right).\] 
Note that a crossing change between bands or a self-crossing change of 
a band can be realized by crossing changes
between $K_0$ and a band as illustrated in Fig. 3.10. Therefore we can
deform each chord into a local
chord by band exchanges and crossing changes between $K_0$ and bands.

When we perform a crossing change between $K_0$ and a $C_p$-band of a
$C_p$-chord ${\cal B}_{p,q}$
with $p\leq k-2$ we introduce a new 
$C_{p+1}$-chord ${\cal B}_{p+1,r}$ and we set $\omega({\cal
B}_{p+1,r})=\omega({\cal B}_{p,q})$ so that the condition $(*)$ still holds
for all subset $P$ of $\{1,...,l\}$. To do this we use Lemma 3.8. 
When we perform a band
exchange between a $C_p$-chord ${\cal B}_{p,q}$ and a $C_r$-chord ${\cal
B}_{r,s}$ with $p+r\leq k-1$ 
we introduce a new
$C_{p+r}$-chord ${\cal B}_{p+r,n}$ and set $\omega({\cal
B}_{p+r,n})=\omega({\cal B}_{p,q})\cup\omega({\cal B}_{r,s})$ so that the
condition $(*)$ still holds for all subset $P$ of
$\{1,...,l\}$. To do this we use Lemma 3.9. Note that the condition 
$\sum_{t\in\omega({\cal B}_{i,j})}k_t\leq i$ still holds for all chords. 

By Lemma 3.8, a crossing change between $K_0$ and a
$C_{k-1}$-band is realized by a
$C_{k}$-move  and therefore does not change the
$C_{k}$-equivalence class.
Similarly, by Lemma 3.9, a band exchange between a 
$C_p$-chord ${\cal B}_{p,q}$ and a $C_r$-chord ${\cal B}_{r,s}$
with $p+r\geq k$ is realized by a $C_{p+r}$-move. 
Therefore, by Remark 2.4, it also does not change the
$C_{k}$-equivalence class.

We note that the process definitely ends at last because a deformation of a
$C_p$-chords does not produces $C_q$-chords for $q\leq p$. $\Box$

\begin{center}
\includegraphics[trim=0mm 0mm 0mm 0mm, width=.7\linewidth]
{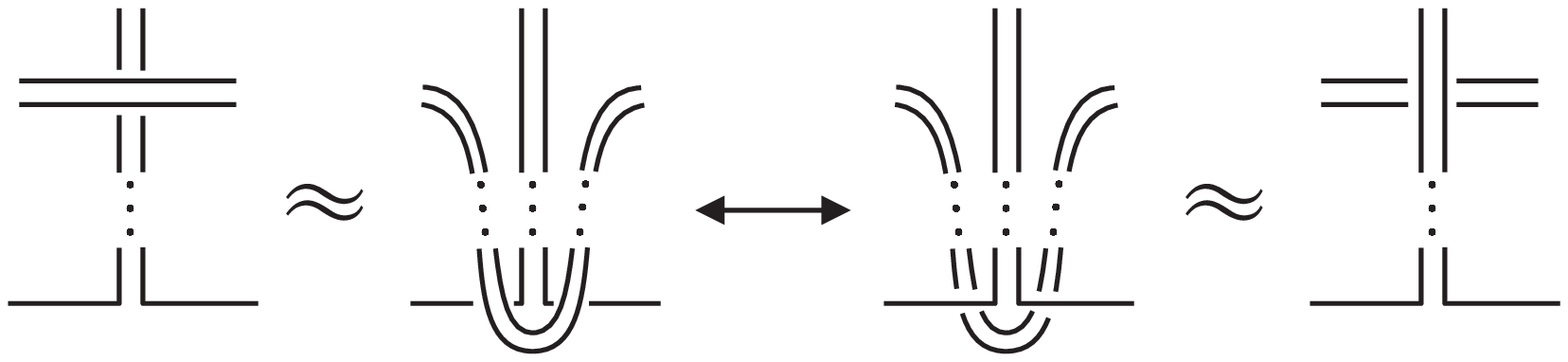}

Fig. 3.10
\end{center}

\medskip
\noindent
{\bf Proof of Theorem 1.4.} It is clear that if $K_1$ and $K_2$ are
$C_k$-equivalent and if $J_1$ and $J_2$ are $C_k$-equivalent then 
the connected sums $K_1\#J_1$ and $K_2\#J_2$ are
$C_k$-equivalent. Thus the binary operation is well-defined. We note that
the ambient isotopy classes
of oriented knots forms a commutative monoid under connected sum of knots
with unit element a trivial
knot. Therefore it is sufficient to show the existence of the inverse
element. Let $K$ be an oriented knot. First we note that
$K$ itself is $C_1$-equivalent to a trivial knot. Suppose that there is an
oriented knot $J$ such that $K\#J$ is $C_{k-1}$-equivalent to a trivial 
knot $O$. Then by Lemma 3.6 we have that
$O$ is a band sum of $K\#J$ and some $C_{k-1}$-link models. We choose a
3-ball $B$ in $S^3$ such that
$B\cap K\#J$ is an unknotted arc in $B$. By an ambient isotopy we deform
the link balls into $B$ and
slide the ends of bands into $B$. Then using Lemma 3.8 we deform $O$ up to
$C_k$-equivalence so that
the whole of the bands are contained in $B$. Then we have that the result
is a connected sum of $K\#J$
and some knot $L$. Namely $K\#J\#L$ is $C_k$-equivalent to $O$. 
Thus $J\#L$ is the desired knot. $\Box$

\medskip
\noindent
{\bf Proof of Theorem 1.5.} 
Let $\xi_k:{\Bbb Z}{\cal K}/({\cal V}(k)+{\cal R}_\#)\rightarrow{\cal
K}/C_k$ be a homomorphism defined by $\xi_k(K)=[K]_{C_k}$ for $K\in{\cal
K}$ where $[K]_{C_k}$ denote
the $C_k$-equivalence class of $K$. 
It follows from Theorem 1.4 that both $\eta_k$ and $\xi_k$ are 
well-defined homomorphisms.
Then it is clear that both $\xi_k\circ\eta_k$ and $\eta_k\circ\xi_k$ are
identities. Therefore both
$\eta_k$ and $\xi_k$ are isomorphisms.
$\Box$

\bigskip
\footnotesize{
 }

\end{document}